\documentclass[a4paper,12pt]{article}
\usepackage{amsthm,amsmath}
\usepackage{amssymb}
\usepackage[all]{xypic}
\usepackage{enumerate}
\usepackage{a4wide}

\swapnumbers
\newtheorem{theorem}{Theorem}[section]
\newtheorem{proposition}[theorem]{Proposition}
\newtheorem{corollary}[theorem]{Corollary}
\newtheorem{lemma}[theorem]{Lemma}
\newtheorem*{theorem*}{Theorem}
\newtheorem*{proposition*}{Proposition}
\newtheorem*{corollary*}{Corollary}
\newtheorem*{lemma*}{Lemma}
\theoremstyle{definition}
\newtheorem{definition}[theorem]{Definition}
\newtheorem{punto}[theorem]{}
\newtheorem{example}[theorem]{Example}

\newtheorem*{remark*}{Remark}

\newtheorem*{definition*}{Definition}

\newcommand{\cat}[1]{\mathcal{#1}}
\newcommand{\coring}[1]{\mathfrak{#1}}
\newcommand{\tensor}[1]{\otimes_{#1}}
\newcommand{\tensfun}[1]{\underset{{#1}}{\otimes}}
\newcommand{\tensdia}[1]{\otimes_{#1}}
\newcommand{\rcomod}[1]{\mathcal{M}^{#1}}
\newcommand{\rmod}[1]{\mathcal{M}_{#1}}
\newcommand{\bcomod}[2]{{}^{#1}\mathcal{M}^{#2}}
\newcommand{\cotensor}[1]{\square_{#1}}
\newcommand{\lcomod}[1]{{}^{#1}\mathcal{M}}

\renewcommand{\hom}[3]{\mathrm{Hom}_{#1}(#2,#3)}

\newcommand{\rend}[2]{\mathrm{End}({#2}_{#1})}
\newcommand{\nat}[2]{\Upsilon_{#1,#2}}
\newcommand{\cohom}[3]{\mathrm{h}_{#1}(#2,#3)}
\newcommand{\adiso}[2]{\Phi_{#1,#2}}

\begin{document}
\title{Separable functors in corings}
\author{J. G\'omez-Torrecillas \\
Departamento de \'{A}lgebra \\
Universidad de Granada\\
E18071 Granada, SPAIN \\
e-mail: \textsf{torrecil@ugr.es}}

\date{\empty}

\maketitle
\begin{abstract}
We develop some basic functorial techniques for the study of the
categories of comodules over corings. In particular, we prove that
the induction functor stemming from every morphism of corings has
a left adjoint, called ad-induction functor. This construction
generalizes the known adjunctions for the categories of Doi-Hopf
modules and entwined modules. The separability of the induction
and ad-induction functors are characterized, extending earlier
results for coalgebra and ring homomorphisms, as well as for
entwining structures.
\end{abstract}
\section*{}
\textit{Keywords and prhases:} coring, comodule, cotensor product
functor, cohom functor, coring homomorphism, induction functor,
ad-induction functor, separable functor \\
\noindent \textit{2000 Mathematics Subject Classification:} 16W30

\section*{Introduction}

The notion of separable functor was introduced by C. N\u ast\u
asescu, M. Van den Bergh and F. Van Oystaeyen in
\cite{Nastasescu/VandenBergh/VanOystaeyen:1989}, where some
applications for group-graded rings were done. This notion fits
satisfactorily to the classical notion of separable algebra over a
commutative ring. Every separable functor between abelian
categories encodes a Maschke's Theorem, which explains the
interest concentrated in this notion within the module-theoretical
developments in recent years. Thus, separable functors have been
investigated in the framework of coalgebras
\cite{Castano/Gomez/Nastasescu:1997}, graded homomorphisms of
rings \cite{DelRio:1992,Castano/Gomez/Nastasescu:1998},
Doi-Koppinen modules
\cite{Caenepeel/alt:1999,Caenepeel/alt:2000unp} or, finally,
entwined modules
\cite{Brzezinski:1999unp,Brzezinski/alt:2000unp1}. These
situations are generalizations of the original study of the
separability for the induction and restriction of scalars functors
associated to a ring homomorphism done in
\cite{Nastasescu/VandenBergh/VanOystaeyen:1989}. It turns out that
that all the aforementioned categories of modules are instances of
comodule categories over suitable corings
\cite{Brzezinski:2000unp}. In fact, the separability of some
fundamental functors relating the category of comodules over a
coring and the underlying category of modules has been studied in
\cite{Brzezinski:2000unp}. Thus, we can expect that the
characterizations obtained in \cite{Brzezinski:1999unp} of the
separability of the induction functor associated to an admissible
morphism of entwining structures and its adjoint generalize to
the corresponding functors stemming from a homomorphism of
corings. This is done in this paper.

To state and prove the separability theorems, I have developed a
basic theory of functors between categories of comodules, making
the arguments independent from the Sweedler's `sigma-notation'. My
plan here is to use purely categorical methods which could be
easily adapted to more general developments of the theory. These
methods had been sketched in \cite{AlTakhman:unp2001} and
\cite{AlTakhman:1999} in the framework of coalgebras over
commutative rings and are expounded in Sections \ref{bicomodulos},
\ref{funtores} and \ref{cohom}. In Section \ref{separable}, I give
a notion of homomorphism of corings, which leads to a pair of
adjoint functors (the induction functor and its adjoint, called
here ad-induction functor). The morphisms of entwining structures
\cite{Brzezinski:1999unp} are instances of homomorphisms of
corings in our setting. Finally, the separability of these
functors is characterized.

We use essentially the categorical terminology of
\cite{Stenstrom:1975}, with the exception of the term $K$--linear
category and functor, for $K$ a commutative ring (see, e.g.
\cite[I.0.2]{SaavedraRivano:1972}). There are, however, some minor
differences: The notation $X \in \cat{A}$ for a category $\cat{A}$
means that $X$ is an object of $\cat{A}$, and the identity
morphism attached to any object $X$ will be represented by the
object itself. The notation $\rmod{K}$ stands for the category of
all unital $K$--modules. The fact that $G$ is a right adjoint to
some functor $F$ will be denoted by $F \dashv G$. For the notion
of separable functor, the reader is referred to
\cite{Nastasescu/VandenBergh/VanOystaeyen:1989}. Finally, let
$f,g : X \rightarrow Y$ be a pair of morphisms of right modules
over a ring $R$, and let $k : K \rightarrow X$ be its equalizer
(that is, the kernel of $f-g$). We will say that a left
$R$--module $Z$ preserves the equalizer of $(f,g)$ if $k
\tensor{R} Z : K \tensor{R} Z \rightarrow X \tensor{R}Z$ is the
equalizer of the pair $(f \tensor{R} Z,g \tensor{R} Z)$. Of
course, every flat module ${}_RZ$ preserves all equalizers.

\section{Bicomodules and the cotensor product
functor}\label{bicomodulos}

I first recall from \cite{Sweedler:1975} the notion of coring. The
concepts of comodule and bicomodule over a coring are
generalizations of the corresponding notions for coalgebras. We
briefly state some basic properties of the cotensor product of
bicomodules. Similar associativity properties were studied in
\cite{Guzman:1989} in the framework of coseparable corings and in
\cite{AlTakhman:unp2001} for coalgebras over commutative rings.

Throughout this paper, $A, A', A'', \dots$ denote associative and
unitary algebras over a commutative ring $K$.

\begin{punto}\label{coringdef}\textbf{Corings.}
 An $A$--\emph{coring} is a three-tuple
$(\coring{C},\Delta_{\coring{C}},\epsilon_{\coring{C}})$
consisting of an $A$--bicomodule $\coring{C}$ and two
$A$--bimodule maps
\begin{equation}
\xymatrix{ \Delta_{\coring{C}} : \coring{C} \ar[r] & \coring{C}
\tensdia{A} \coring{C} & \epsilon_{\coring{C}} : \coring{C} \ar[r]
& A}
\end{equation}
such that the diagrams
\begin{equation}
\xymatrix{ \coring{C} \ar^{\Delta_{\coring{C}}}[rr]
\ar_{\Delta_{\coring{C}}}[d]  & & \coring{C} \tensdia{A}
\coring{C}
\ar^{\coring{C} \tensfun{A} \Delta_{\coring{C}}}[d] \\
\coring{C} \tensdia{A} \coring{C} \ar^{\Delta_{\coring{C}}
\tensfun{A} \coring{C}}[rr] & & \coring{C} \tensdia{A} \coring{C}
\tensdia{A} \coring{C}}
\end{equation}
and
\begin{equation}
\xymatrix{ \coring{C} \ar^{\Delta_{\coring{C}}}[r] \ar_{\cong}[dr]
& \coring{C} \tensdia{A} \coring{C} \ar^{\coring{C} \tensfun{A}
\epsilon_{\coring{C}}}[d] \\
& \coring{C} \tensdia{A} A } \qquad \xymatrix{ \coring{C}
\ar^{\Delta_{\coring{C}}}[r] \ar_{\cong}[dr] & \coring{C}
\tensdia{A} \coring{C}
\ar^{\epsilon_{\coring{C}} \tensfun{A} \coring{C}}[d] \\
& A \tensdia{A} \coring{C} }
\end{equation}
commute.
\end{punto}

From now on, $\coring{C}, \coring{C}', \coring{C}'', \dots$ will
denote corings over $A, A', A'', \dots$, respectively.

\begin{punto}\label{comoddef}\textbf{Comodules.}
A \emph{right $\coring{C}$--comodule} is a pair $(M, \rho_M)$
consisting of a right $A$--module $M$ and an $A$--linear map
$\rho_M : M \rightarrow M \tensor{A} \coring{C}$ such that the
diagrams
\begin{equation}
\xymatrix{ M \ar^{\rho_M}[rr] \ar_{\rho_M}[d] & {} & M \tensdia{A}
\coring{C} \ar^{M \tensfun{A} \Delta_{\coring{C}}}[d] \\
M \tensdia{A} \coring{C} \ar^{\rho_M \tensfun{A} \coring{C}}[rr] &
{} & M \tensdia{A} \coring{C} \tensdia{A} \coring{C} } \qquad
\xymatrix{ M \ar^{\rho_M}[r] \ar_{\cong}[dr] & M \tensdia{A}
\coring{C}
\ar^{M \tensfun{A} \epsilon_{\coring{C}}}[d] \\
& M \tensdia{A} A }
\end{equation}
commute. Left $\coring{C}$--comodules are similarly defined; we
use the notation $\lambda_M$ for their structure maps. A
\emph{morphism} of right $\coring{C}$--comodules $(M,\rho_M)$ and
$(N, \rho_N)$ is an $A$--linear map $f : M \rightarrow N$ such
that the following diagram is commutative.
\begin{equation}
\xymatrix{
M \ar^f[rr] \ar^{\rho_M}[d] & & N \ar^{\rho_N}[d] \\
M \tensdia{A} \coring{C} \ar^{f \tensfun{A} \coring{C}}[rr] & & N
\tensdia{A} \coring{C}}
\end{equation}
The $K$--module of all right $\coring{C}$--comodule morphisms from
$M$ to $N$ is denoted by $\hom{\coring{C}}{M}{N}$. The $K$--linear
category of all right $\coring{C}$--comodules will be denoted by
$\rcomod{\coring{C}}$. When $\coring{C} = A$, with the trivial
coring structure, the category $\rcomod{A}$ is just the category
of all right $A$--modules, which is `traditionally' denoted by
$\rmod{A}$.

Coproducts and cokernels in $\rcomod{\coring{C}}$ do exist and
can be already computed in $\rcomod{A}$. Therefore,
$\rcomod{\coring{C}}$ has arbitrary inductive limits. If
${}_A\coring{C}$ is a flat module, then $\rcomod{\coring{C}}$ is
easily proved to be an abelian category.
\end{punto}

\begin{punto}\textbf{Bicomodules.}
Let $\rho_M : M \rightarrow M \tensor{A} \coring{C}$ be a comodule
structure over an $A'-A$--bimodule $M$, and assume that $\rho_M$
is $A'$--linear. For any right $A'$--module $X$, the right
$A$--linear map $X \tensfun{A'} \rho_{M} : X \tensor{A'} M
\rightarrow X \tensor{A'} M \tensor{A} \coring{C}$ makes $X
\tensor{A'} M$ a right $\coring{C}$--comodule. This leads to an
additive functor $- \tensor{A'} M : \rmod{A'} \rightarrow
\rcomod{\coring{C}}$. When $A' = A$ and $M = \coring{C}$, the
functor $- \tensor{A} \coring{C}$ is left adjoint to the
forgetful functor $U_A : \rcomod{\coring{C}} \rightarrow
\rcomod{A}$ (see \cite[Proposition 3.1]{Guzman:1989}, \cite[Lemma
3.1]{Brzezinski:2000unp}). Now assume that the $A'-A$--bimodule
$M$ is also a left $\coring{C}'$--comodule with structure map
$\lambda_M : M \rightarrow \coring{C}' \tensor{A} M$. It is clear
that $\rho_M : M \rightarrow M \tensor{A} \coring{C}$ is a
morphism of left $\coring{C}'$--comodules if and only if
$\lambda_M : M \rightarrow \coring{C}' \tensor{A'} M$ is a
morphism of right $\coring{C}$--comodules. In this case, we say
that $M$ is a $\coring{C}'-\coring{C}$--bicomodule. The
$\coring{C}'-\coring{C}$--bicomodules are the objects of a
$K$--linear category $\bcomod{\coring{C}'}{\coring{C}}$ whose
morphisms are those $A'-A$--bimodule maps which are morphisms of
$\coring{C}'$--comodules and of $\coring{C}$--comodules. Some
particular cases are now of interest. For instance, when
$\coring{C}' = A'$, the objects of the category
$\bcomod{A'}{\coring{C}}$ are the $A'-A$--bimodules with a right
$\coring{C}$--comodule structure $\rho_M : M \rightarrow M
\tensor{A} \coring{C}$ which is $A'$--linear.
\end{punto}

\begin{punto}\label{cotensorfun}\textbf{The cotensor product.}
Let $M \in \bcomod{\coring{C}'}{\coring{C}}$ and $N \in
\bcomod{\coring{C}}{\coring{C}''}$. We consider $M \tensor{A} N$
and $M \tensor{A} \coring{C} \tensor{A} N$  as
$\coring{C}'-\coring{C}''$--bicomodules with structure maps
\begin{equation}
\xymatrix{ M \tensdia{A} N \ar^{M \tensfun{A} \rho_N}[rr] & & M
\tensdia{A} N \tensdia{A''} \coring{C}'' & &  M \tensdia{A} N
\ar^{\lambda_M \tensfun{A} M }[rr] & & \coring{C}' \tensdia{A'} M
\tensdia{A} N}
\end{equation}
\begin{equation}
\xymatrix{ M \tensdia{A} \coring{C} \tensdia{A} N \ar^{M
\tensfun{A} \coring{C} \tensfun{A} \rho_N}[rrr] & & & M
\tensdia{A} \coring{C} \tensdia{A} N \tensdia{A''} \coring{C}''}
\end{equation}
\begin{equation}
\xymatrix{
 M \tensdia{A} \coring{C} \tensdia{A} N \ar^{\lambda_M \tensfun{A}
\coring{C} \tensfun{A} N}[rrr] & & & \coring{C}' \tensdia{A'} M
\tensdia{A} \coring{C} \tensdia{A} N}
\end{equation}
The map
\begin{equation}\label{precotensor}
\xymatrix{ M \tensdia{A} N \ar^{\rho_M \tensfun{A} N - M
\tensfun{A} \lambda_N}[rrr] & & & M \tensdia{A} \coring{C}
\tensdia{A} N}
\end{equation}
is then a $\coring{C}'-\coring{C}''$-bicomodule map. Let $M
\cotensor{\coring{C}} N$ denote the kernel of \eqref{precotensor}.
 If $\coring{C}'_{A'}$ and ${}_{A''}\coring{C}''$ preserve the equalizer of
 $(\rho_M \tensor{A} N,M \tensor{A}
\lambda_N)$, then $M \cotensor{\coring{C}} N$ is both a
$\coring{C}'$ and a $\coring{C}''$--subcomodule of $M \tensor{A}
N$ and, hence, it is a $\coring{C}'-\coring{C}''$--bicomodule.
\begin{proposition*}
Assume that $\coring{C}'_{A'}$ and ${}_{A''}\coring{C}''$
preserve the equalizer of $(\rho_M \tensor{A} N,M \tensor{A}
\lambda_N)$ for every $M \in \bcomod{\coring{C}'}{\coring{C}}$
and $N \in \bcomod{\coring{C}}{\coring{C}''}$. We have an
additive bifunctor
\begin{equation}\label{cotenbifun}
\xymatrix{
 - \cotensor{\coring{C}} - :
\bcomod{\coring{C}'}{\coring{C}} \times
\bcomod{\coring{C}}{\coring{C}''} \ar[rr] & &
\bcomod{\coring{C}'}{\coring{C}''}}
\end{equation}
In particular, the cotensor product bifunctor \eqref{cotenbifun}
is defined when $\coring{C}'_{A'}$ and ${}_{A''}\coring{C}''$ are
flat modules or when $\coring{C}$ is a coseparable $A$--coring in
the sense of \cite{Guzman:1989}.
\end{proposition*}

 In the special case $\coring{C}' = A'$ and $\coring{C}'' =
A''$, we have the bifunctor
\begin{equation}
\xymatrix{
- \cotensor{\coring{C}} - : \bcomod{A'}{\coring{C}}
\times \bcomod{\coring{C}}{A''} \ar[rr] & & \bcomod{A'}{A''}}
\end{equation}
and, if we further assume $A' = A'' = K$, we have the bifunctor
\begin{equation}
\xymatrix{ - \cotensor{\coring{C}} - : \rcomod{\coring{C}} \times
\lcomod{\coring{C}} \ar[rr] & & \rmod{K}}.
\end{equation}
\end{punto}

\begin{punto}\label{coass1}\textbf{Compatibility between tensor
and cotensor.} Let $M \in \bcomod{\coring{C}'}{\coring{C}}$ and $N
\in \bcomod{\coring{C}}{\coring{C}''}$ be bicomodules. For any
right $A'$--module $W$, consider the commutative diagram
\begin{equation}\label{primeaso}
\xymatrix{ & W \tensdia{A'} (M \cotensor{\coring{C}} N) \ar[r]
\ar^{\psi}[d] & W \tensdia{A'} (M \tensdia{A} N) \ar[r]
\ar^{\cong}[d] & W
\tensdia{A'} (M \tensdia{A} \coring{C} \tensdia{A} N) \ar^{\cong}[d] \\
0 \ar[r] & (W \tensdia{A'} M) \cotensor{\coring{C}} N \ar[r] & (W
\tensdia{A'} M) \tensdia{A} N \ar[r] & (W \tensdia{A'} M)
\tensdia{A} \coring{C} \tensdia{A} N },
\end{equation}
where $\psi$ is given by the universal property of the kernel in
the second row. This leads to the following
\begin{lemma*}It follows from \eqref{primeaso} that
$W_{A'}$ preserves the equalizer of $(\rho_M \tensor{A} N, M
\tensor{A} \lambda_N)$ if and only if $\psi : W \tensor{A'} ( M
\cotensor{\coring{C}} N) \cong (W \tensor{A'} M)
\cotensor{\coring{C}} N$. In particular, $\psi$ is an isomorphism
if $W_{A'}$ is flat.
\end{lemma*}

Next, we prove a basic fact concerning with the associativity of
cotensor product.
\begin{proposition*}
Let $M \in \bcomod{\coring{C}'}{\coring{C}}, N \in
\bcomod{\coring{C}}{\coring{C}''}, L \in
\bcomod{\coring{C}'''}{\coring{C}'}$.
Assume that $\coring{C}'_{A'}, L_{A'}, L \tensor{A'}
\coring{C}'_{A'}$ and ${}_{A''}\coring{C}''$ preserve the
equalizer of $(\rho_M \tensor{A} N, M \tensor{A} \lambda_N)$, and
that ${}_{A}\coring{C}, {}_AN, {}_A\coring{C} \tensor{A} N$ and
$\coring{C}'''_{A'''}$ preserve the equalizer of $(\rho_L
\tensor{A'} M, L \tensor{A'} \lambda_M)$. Then we have a canonical
isomorphism of $\coring{C}'''-\coring{C}''$--bicomodules
\[
L \cotensor{\coring{C}'} (M \cotensor{\coring{C}} N) \cong (L
\cotensor{\coring{C'}} M) \cotensor{\coring{C}} N
\]
\end{proposition*}
\begin{proof}
Since $\coring{C}'_{A'}$ and ${}_{A''}\coring{C}''$ preserve the
equalizer of $(\rho_M \tensor{A} N, M \tensor{A} \lambda_N)$, we
know that it is a $\coring{C}'-\coring{C}''$--subbicomodule.
Analogously, $L \cotensor{\coring{C}'} M$ is a
$\coring{C}'''-\coring{C}$--subbicomodule of $L \tensor{A'} M$ In
the commutative diagram
\[
\xymatrix{ 0 \ar[r] & (L \cotensor{\coring{C}'} M)
\cotensor{\coring{C}} N \ar[r] \ar[d] & (L \tensor{A'} M)
\cotensor{\coring{C}} N \ar[r] \ar[d] & (L
\tensor{A'} \coring{C}' \tensor{A'} M) \cotensor{\coring{C}} N \ar[d] \\
0 \ar[r] & (L \cotensor{\coring{C}'} M) \tensor{A} N \ar[r] & (L
\tensor{A'} M ) \tensor{A} N \ar[r] & (L \tensor{A'} \coring{C}'
\tensor{A'} M) \tensor{A} N  }
\]
the second row is exact because ${}_AN$ preserves the equalizer
of $(\rho_L \tensor{A'} M, L \tensor{A'} \lambda_M)$. The
exactness of the first row is then deduced by using that
${}_A\coring{C} \tensor{A} N$ is assumed to preserve the
equalizer of $(\rho_L \tensor{A'} M, L \tensor{A'} \lambda_M)$.
Now, consider the commutative diagram with exact rows:
\[
\xymatrix{0 \ar[r] & (L \cotensor{\coring{C}'} M)
\cotensor{\coring{C}} N \ar[r] & (L \tensor{A'} M)
\cotensor{\coring{C}} N \ar[r] & (L
\tensor{A'} \coring{C}' \tensor{A'} M) \cotensor{\coring{C}} N \\
0 \ar[r] & L \cotensor{\coring{C}'} ( M \cotensor{\coring{C}} N)
\ar[r] \ar^{\psi_1}[u] & L \tensor{A'} (M \cotensor{\coring{C}}
N) \ar[r] \ar^{\psi_2}[u] & L \tensor{A'} \coring{C}' \tensor{A'}
(M \cotensor{\coring{C}} N) \ar^{\psi_3}[u]}
\]
Lemma \ref{coass1} gives the isomorphisms $\psi_2$ and $\psi_3$,
which induce the isomorphism $\psi_1$.
\end{proof}
\end{punto}

\section{Functors between comodule categories}\label{funtores}

This section contains technical facts concerning with $K$--linear
functors between categories of comodules over corings. Part of
these tools were first developed for coalgebras over commutative
rings in \cite{AlTakhman:unp2001} and \cite{AlTakhman:1999}.
Roughly speaking, I prove an analogue to Watts theorem, which
allow to represent good enough functors as cotensor product
functors. I also include a result which states that, under mild
conditions, a natural transformation gives a bicomodule morphism
at any bicomodule. This will be used in the statement and proof
of our separability theorems in Section \ref{separable}.

Let $\coring{C}$, $\coring{D}$ be corings over $K$--algebras $A$
and $B$, respectively, and consider a $K$--linear functor
\[
\xymatrix{ F : \rcomod{\coring{C}} \ar[r] & \rcomod{\coring{D}}}
\]
\begin{punto}\label{nat1}
Let $T$ be a $K$--algebra. For every $M \in
\bcomod{T}{\coring{C}}$, consider the homomorphism of
$K$--algebras
\begin{equation}\label{homobasico}
 T \cong \rend{T}{T} \rightarrow  \hom{\coring{C}}{T \tensdia{T}
M}{T \tensdia{T} M} \cong \hom{\coring{C}}{M}{M} \rightarrow
\hom{\coring{D}}{F(M)}{F(M)}
\end{equation}
which induces a left $T$--module structure over $F(M)$ such that
$F(M)$ becomes a $T-\coring{D}$--bicomodule. We have two
$K$--linear functors
\[
\xymatrix{ - \tensdia{T} F(-), F(- \tensdia{T} -) : \rcomod{T}
\times \bcomod{T}{\coring{C}} \ar[r] & \bcomod{}{\coring{D}} }
\]
We shall construct a natural transformation
\[
\xymatrix{ \nat{-}{-} : - \tensdia{T} F(-) \ar[r] & F(-
\tensdia{T} -) }
\]
Let $\nat{T}{M}$ be the unique isomorphism of
$\coring{D}$--comodules making commutative the diagram
\[
\xymatrix{ T \tensdia{T} F(M) \ar^{\nat{T}{M}}[rr] \ar^{\cong}[dr]
& &
F(T \tensdia{T} M) \ar^{\cong}[dl] \\
& F(M)  & }
\]
To prove that $\nat{T}{M}$ is natural at $T$, consider a
homomorphism $f : T \rightarrow T$ of right $T$--modules and
define $g : F(M) \rightarrow F(M)$ by $g(x) = f(1)x$ for every $x
\in F(M)$. Since $g$ is just the image under \eqref{homobasico} of
$f(1) \in T$, it follows that $g$ is a morphism of right
$\coring{D}$--comodules. Moreover, $g$ makes the following
diagram commutative:
\[
\xymatrix{
F(M) \ar^{g}[rr] & & F(M) \\
F(T \tensdia{T} M) \ar^{\cong}[u] \ar^{F(f \tensfun{T} M)}[rr] & &
F(T \tensdia{T} M) \ar^{\cong}[u]}
\]
In the diagram
\[
\xymatrix{
& F(M) \ar^{g}[rr]  & & F(M) \\
T \tensdia{T} F(M) \ar^{\cong}[ur] \ar^{f \tensfun{T} F(M)}[rr]
\ar_{\nat{T}{M}}[dd] &
& T \tensdia{T} F(M) \ar^{\cong}[ur] \ar_{\nat{T}{M}}[dd] & \\
 &  & &  \\
 F(T \tensdia{T} M) \ar_{\cong}[uuur] \ar^{F(f \tensfun{T} M)}[rr] & &
 F(T \tensdia{T} M) \ar_{\cong}[uuur]}
\]
the commutativity of the front rectangle, which gives that
$\nat{T}{M}$ is natural, follows from the commutativity of the
rest of the diagram. From Mitchell's theorem (\cite[Theorem 3.6.5]
{Popescu:1973}) we obtain a natural transformation
\[
\xymatrix{ \nat{-}{M} : - \tensdia{T} F(M) \ar[r] & F( -
\tensdia{T} M)}.
\]
Moreover, if $F$ preserves coproducts (resp. direct limits, resp.
inductive limits) then $\nat{X}{M}$ is an isomorphism for $X_T$
projective (resp. flat, resp. any right $T$-module).
\end{punto}

\begin{proposition}\label{natpreserves} If the functor $F$
preserves coproducts, then $\nat{-}{-} : - \tensor{T} F(-)
\rightarrow F(- \tensor{T} -)$ is a natural transformation.
Moreover, if $F$ preserves direct limits, then $\nat{X}{-}$ is a
natural isomorphism for every flat right $T$--module $X$.
Finally, if $F$ preserves inductive limits, then $\nat{-}{-}$ is
a natural isomorphism.
\end{proposition}
\begin{proof}
By \ref{nat1}, $\nat{-}{M}$ is natural for every $M \in
\bcomod{T}{\coring{C}}$. Thus, we have only to show that
$\nat{X}{-}$ is natural for every $X \in \rcomod{T}$. We argue
first for $X = T$. Let $f : M \rightarrow N$ be a homomorphism in
$\bcomod{T}{\coring{C}}$. From the diagram
\[
\xymatrix{
 & F(M) \ar^{F(f)}[rr]  & & F(N) \\
T \tensdia{T} F(M) \ar^{\cong}[ur] \ar^{T \tensfun{T} F(f)}[rr]
\ar_{\nat{T}{M}}[dd] &
& T \tensdia{T} F(M) \ar^{\cong}[ur] \ar_{\nat{T}{N}}[dd] & \\
 &  & &  \\
 F(T \tensdia{T} M) \ar_{\cong}[uuur] \ar^{F(T \tensdia{T} f)}[rr] & &
 F(T \tensdia{T} N) \ar_{\cong}[uuur]}
\]
we get that $\nat{A}{-}$ is natural. Now, use a free presentation
$T^{(\Omega)} \rightarrow X$ to obtain that $\nat{X}{-}$ is
natural for a general $X_T$. The rest of the statements are easily
derived from this.
\end{proof}

\begin{lemma}\label{FunoFdos}
Let $\eta : F_1 \rightarrow F_2$ a natural transformation, where
$F_1, F_2 : \rcomod{\coring{C}} \rightarrow \rcomod{\coring{D}}$
are $K$--linear functors which preserve coproducts.
\begin{enumerate}
\item
For every $M \in \bcomod{T}{\coring{C}}$, $\eta_M : F_1(M)
\rightarrow F_2(M)$ is a $T-\coring{D}$--bicomodule homomorphism.
\item Given $X \in \rcomod{T}$ and $M \in \bcomod{T}{\coring{C}}$, the diagram
\begin{equation}\label{etas}
\xymatrix{ F_1(X \tensdia{T} M) \ar^{\eta_{X \tensdia{T} M}}[rr] &
&
F_2(X \tensdia{T} M) \\
X \tensdia{T} F_1(M) \ar^{X \tensfun{T} \eta_M}[rr]
\ar^{\nat{X}{M}}[u]& & X \tensdia{T} F_2(M) \ar^{\nat{X}{M}}[u]}
\end{equation}
is commutative.
\end{enumerate}
\end{lemma}
\begin{proof}
We need just to prove that \eqref{etas} commutes for $X = T$. In
this case, the diagram can be factored out as
\[
\xymatrix{ F_1(T \tensdia{T} M) \ar^{\eta_{T \tensdia{T} M}}[rr]
\ar^{\cong}[dr]& & F_2(T \tensdia{T} M) \ar^{\cong}[dr]&
\\
 & F_1(M) \ar^(.3){\eta_M}[rr] & \ar[u] & F_2(M)
\\
T \tensdia{T} F_1(M) \ar^{\cong}[ur] \ar^{T \tensfun{T}
\eta_M}[rr] \ar^{\nat{T}{M}}[uu] & & T \tensdia{T} F_2(M)
\ar^{\cong}[ur] \ar@{-}^(.6){\nat{T}{M}}[u]& }
\]
Since all trapezia and triangles commute, the back rectangle does,
as desired.
\end{proof}

\begin{lemma}\label{formula}
Let $T, S$ be $K$--algebras and assume that $F :
\rcomod{\coring{C}} \rightarrow \rcomod{\coring{D}}$ preserves
coproducts. Given $X \in \rcomod{S}, Y \in \bcomod{S}{T}$ and $M
\in \bcomod{T}{\coring{C}}$, the following formula holds
\[
\nat{X}{Y \tensdia{T} M} \circ (X \tensfun{S} \nat{Y}{M}) = \nat{X
\tensor{S} Y}{M}
\]
\end{lemma}
\begin{proof}
The equality will be first proved for $X = S$. Consider the
diagram
\[
\xymatrix{ S \tensdia{S} F(Y \tensdia{T} M) \ar^{\cong}[rr]
\ar^{\nat{S}{Y \tensdia{T} M}}[dd] & & F(Y
\tensdia{T} M) \ar@{=}'[d][dd] & \\
 & S \tensdia{S} Y \tensdia{T} F(M) \ar^{S \tensfun{S} \nat{Y}{M}}[ul]
 \ar^{\nat{S \tensdia{S} Y}{M}}[dl] \ar^{\cong}[rr] & & Y \tensdia{T} F(M)
 \ar^{\nat{Y}{M}}[ul] \ar^{\nat{Y}{M}}[dl] \\
 F(S \tensdia{S} Y \tensdia{T} M) \ar^{\cong}[rr] & & F(Y \tensdia{T} M) & }
\]
The back rectangle is commutative by definition of $\nat{S}{Y
\tensor{T} M}$, while the other two parallelograms are commutative
because $\nat{-}{-}$ is natural. Therefore, the right triangle is
commutative. The equality is now easily extended for $X =
S^{(\Omega)}$ and, by using a free presentation $S^{(\Omega)}
\rightarrow X \rightarrow 0$, for any $X$.
\end{proof}

\begin{punto}\textbf{Natural transformations and
bicomodule morphisms.}\label{hacebicom}  Let $M \in
\bcomod{\coring{C}'}{\coring{C}}$ a bicomodule. A functor $F:
\rcomod{\coring{C}} \rightarrow \rcomod{\coring{D}}$ is said to
be \emph{$M$--compatible} if $\nat{\coring{C}'}{M}$ and
$\nat{\coring{C}'\tensor{A'}\coring{C}'}{M}$ are isomorphisms. By
Proposition \ref{natpreserves}, the functor $F$ is
$M$--compatible for every bicomodule $M$ either if the functor
$F$ preserves inductive limits or $\coring{C}_A'$ is flat and $F$
preserves direct limits. In case that $F$ is $M$--compatible,
define $\lambda_{F(M)}$ as the unique $A'$--linear map making
commutative the diagram
\[
\xymatrix{
F(M) \ar^{\lambda_{F(M)}}[rr] \ar_{F(\lambda_M)}[dr] & &
\coring{C}' \tensdia{A'} F(M) \ar^{\nat{\coring{C}'}{M}}[dl] \\
& F(\coring{C}' \tensdia{A'} M) & }
\]
\begin{proposition*}
Let $F$ be an $M$--compatible functor which preserves coproducts.
The $A'$--linear map $\lambda_{F(M)}$ is a left
$\coring{C}'$--comodule structure on $F(M)$ such that $F(M)$
becomes a $\coring{C}'-\coring{D}$--bicomodule. Moreover, given
$F_1, F_2 : \rcomod{\coring{C}} \rightarrow \rcomod{\coring{D}}$
$M$--compatible functors and a natural transformation $\eta: F_1
\rightarrow F_2$, the map $\eta_M : F_1(M) \rightarrow F_2(M)$ is
a $\coring{C}'-\coring{D}$--bicomodule homomorphism.
\end{proposition*}
\begin{proof}
In order to prove that the coaction $\lambda_{F(M)}$ is
coassociative, let us consider the diagram:
\[
\xymatrix{ &  F(M) \ar[dl]_{\lambda_{F(M)}}
\ar[rr]^{\lambda_{F(M)}} \ar@{=}[d] & & \coring{C}' \tensor{A'}
F(M) \ar[dl]|{\Delta_{\coring{C}'} \tensfun{A'} F(M)}
\ar[ddd]^{\Upsilon_{\coring{C}',M}} \\
\coring{C}' \tensor{A'} F(M) \ar[rr]|(.6){\coring{C}' \tensfun{A'}
\lambda_{F(M)}} \ar[dr]|{\coring{C}' \tensfun{A'} F(M)}
\ar[ddd]^{\Upsilon_{\coring{C}',M}} & \ar@{=}[d] & \coring{C} '
\tensor{A'} \coring{C}' \tensor{A'} F(M) \ar[dl]|{\coring{C}'
\tensfun{A'} \Upsilon_{\coring{C}',M}}
\ar[ddd]^{\Upsilon_{\coring{C}' \tensor{A'} \coring{C}',M}} & \\
 &  \coring{C}' \tensor{A'} F(\coring{C}' \tensor{A'} M)
 \ar[ddr]|(.35){\Upsilon_{\coring{C}',\coring{C}' \tensor{A'} M}} \ar@{=}[d] & & \\
 & F(M) \ar[dl]^{F(\lambda_{M})} \ar[rr]^(.4){F(\lambda_{M})}  & & F(\coring{C}' \tensor{A'} M)
 \ar[dl]^{F(\Delta_{\coring{C}'} \tensfun{A'} M)} \\
F(\coring{C}' \tensor{A'} M) \ar[rr]^{F(\coring{C'} \tensfun{A'}
\lambda_M)} & & F(\coring{C}' \tensor{A'} \coring{C}' \tensor{A}'
M) & }
\]
We want to see that the top side is commutative. Since $F$ is
assumed to be $M$--compatible, we have just to prove that the
mentioned side is commutative after composing with the isomorphism
$\nat{\coring{C}' \tensor{A'} \coring{C}'}{M}$. This is deduced
by using Lemma \ref{formula}, in conjunction with the naturality
of $\nat{-}{M}$ and the very definition of $\lambda_{F(M)}$. The
counitary property is deduced from the commutative diagram
\[
\xymatrix{ A' \tensdia{A'} F(M) \ar^{\nat{A'}{M}}[rr] & & F(A'
\tensdia{A'} M)\\
& F(M) \ar_{\cong}[ul] \ar^{\cong}[ur] \ar_{\lambda_{F(M)}}[dl] \ar^{F(\lambda_M)}[dr]&  \\
\coring{C}' \tensdia{A'} F(M) \ar^{\epsilon_{\coring{C}'}
\tensfun{A'} F(M)}[uu] \ar^{\nat{\coring{C}'}{M}}[rr] & &
F(\coring{C}' \tensdia{A'} M) \ar_{F(\epsilon_{\coring{C}'}
\tensfun{A'} M)}[uu] }
\]
To prove the second statement, let us consider the diagram
\[
\xymatrix{ F_1(M) \ar^{\eta_M}[rr] \ar^{\lambda_{F_1(M)}}[dd]
\ar^{F_1(\lambda_M)}[dr]& & F_2(M) \ar^{F_2(\lambda_M)}[dr]
\ar@{-}^(.7){\lambda_{F_2(M)}}[d] & \\
& F_1(\coring{C}' \tensdia{A'} M) \ar^(.3){\eta_{\coring{C}'
\tensfun{A'} M}}[rr]
\ar^{\nat{\coring{C}'}{M}}[dl]& \ar[d] & F_2(\coring{C}' \tensdia{A} M) \ar^{\nat{\coring{C}'}{M}}[dl]\\
 \coring{C}' \tensdia{A'} F_1(M) \ar^{\coring{C}' \tensfun{A'}
\eta_M}[rr]& & \coring{C}' \tensdia{A'} F_2(M) &}
\]
Both triangles commute by definition of $\lambda_{F_1(M)}$ and
$\lambda_{F_2(M)}$, and the upper trapezium is commutative because
$\eta$ is natural. The bottom trapezium commutes by Proposition
\ref{FunoFdos}. Therefore, the back rectangle is commutative,
which just says that $\eta_M$ is a morphism of left
$\coring{C}'$--comodules. This finishes the proof.
\end{proof}
\end{punto}

\begin{theorem}
Assume that $\coring{C}_A$ is flat. If $F : \rcomod{\coring{C}}
\rightarrow \rcomod{\coring{D}}$ is exact and preserves direct
limits (e.g. if $F$ is an equivalence of categories), then $F$ is
naturally isomorphic to $- \cotensor{\coring{C}} F(\coring{C})$.
\end{theorem}
\begin{proof}
Let $\rho_N \rightarrow N \tensor{A} \coring{C}$ be a right
$\coring{C}$--comodule. We have the following diagram with exact
rows
\[
\xymatrix{ 0 \ar[r] & N \cotensor{\coring{C}} F(\coring{C})
\ar^{\cong}@{-->}[d]\ar[r] & N \tensdia{A} F(\coring{C})
\ar^{\rho_N \tensfun{A} F(\coring{C}) - N \tensfun{A}
\lambda_{F(\coring{C})}}[rrrr]
\ar_{\cong}^{\nat{N}{\coring{C}}}[d]& & & & N \tensdia{A}
\coring{C} \tensdia{A} F(\coring{C})
\ar_{\cong}^{\nat{N \tensdia{A} A}{\coring{C}}}[d] \\
0 \ar[r] & F(N) \ar[r] & F(N \tensdia{A} \coring{C}) \ar^{F(\rho_N
\tensfun{A} \coring{C} - N \tensfun{A} \Delta_{\coring{C}}
)}[rrrr] & & & & F(N \tensdia{A} \coring{C} \tensdia{A}
\coring{C})}
\]
where the desired isomorphism is given by the universal property
of the kernel.
\end{proof}

\section{Co-hom functors}\label{cohom}

This section contains a quick study of the left adjoint to a
cotensor product functor, if it does exist. The presentation is
inspired from the one given in \cite[1.8]{Takeuchi:1977} for
coalgebras over a field.

Let $\coring{C}$, $\coring{D}$ be corings over $K$--algebras $A$
and $B$, respectively.

\begin{definition}
A bicomodule $N \in \bcomod{\coring{C}}{\coring{D}}$ is said to be
\emph{quasi-finite} as a right $\coring{D}$--comodule if the
functor $- \tensor{A} N : \rcomod{A} \rightarrow
\rcomod{\coring{D}}$ has a left adjoint $\cohom{\coring{D}}{N}{-}
: \rcomod{\coring{D}} \rightarrow \rcomod{A}$. This functor is
called the \emph{co-hom} functor associated to $N$.
\end{definition}

The natural isomorphism which gives the adjunction will be denoted
by
\begin{equation}\label{adjuncion}
\xymatrix{ \adiso{X}{Y} : \hom{A}{\cohom{\coring{D}}{N}{X}}{Y}
\ar[rr] & & \hom{\coring{D}}{X}{Y \tensdia{A} N},}
\end{equation}
for $Y \in \rcomod{A}$, $X \in \rcomod{\coring{D}}$.

\begin{punto}\label{casichico}
Let $\theta : Id \rightarrow \cohom{\coring{D}}{N}{-} \tensor{A}
N$ be the unit of the adjunction \eqref{adjuncion}. Therefore,
the isomorphism $\adiso{X}{Y}$ is given by the assignment $f
\mapsto (f \tensor{A} N) \theta_X$. In particular, the map
\[
\xymatrix{ X \ar^(.3){\theta_X}[r]  & \cohom{\coring{D}}{N}{X}
\tensdia{A} N \ar^{id \tensfun{A} \lambda_N }[rr]  &
&\cohom{\coring{D}}{N}{X} \tensdia{A} \coring{C} \tensdia{A} N}
\]
determines an $A$--linear map
\[
\xymatrix{ \rho_{\cohom{\coring{D}}{N}{X}} :
\cohom{\coring{D}}{N}{X} \ar[r] & \cohom{\coring{D}}{N}{X}
\tensdia{A} \coring{C}}
\]
such that $(id \tensor{A} \lambda_N) \theta_X =
(\rho_{\cohom{\coring{D}}{N}{X}} \tensor{A} id) \theta_X$. The
coaction $\rho_{\cohom{\coring{D}}{N}{X}}$ makes
$\cohom{\coring{D}}{N}{X}$ a right $\coring{C}$--comodule.
Therefore we have a functor $\cohom{\coring{D}}{N}{-} :
\rcomod{\coring{D}} \rightarrow \rcomod{\coring{C}}$.
\end{punto}

\begin{punto}\label{adjcohom}
The following is a basic tool in our investigation.

\begin{proposition*}
Let $N$ be a $\coring{C}-\coring{D}$--bicomodule. Assume that
${}_B\coring{D}$ preserves the equalizer of $(\rho_Y \tensor{A} N,
Y \tensor{A} \lambda_N)$ for every right $\coring{C}$--comodule
$Y$ (e.g., ${}_B\coring{D}$ is flat or $\coring{C}$ is a
coseparable $A$--coring in the sense of \cite{Guzman:1989}).
\begin{enumerate}
\item
If $N$ is quasi-finite as a right $\coring{D}$--comodule, then the
natural isomorphism \eqref{adjuncion} restricts to an isomorphism
$\hom{\coring{C}}{\cohom{\coring{D}}{N}{X}}{Y} \cong
\hom{\coring{D}}{X}{Y \cotensor{\coring{C}} N}$. Therefore,
$\cohom{\coring{D}}{N}{-}$ is left adjoint to $-
\cotensor{\coring{C}} N$.
\item Conversely, if
$- \cotensor{\coring{C}} N : \rcomod{\coring{C}} \rightarrow
\rcomod{\coring{D}}$ has a left adjoint, then $N$ is quasi-finite
as a right $\coring{C}$--comodule.
\end{enumerate}
\end{proposition*}
\begin{proof}
1. We need to prove that if $f \in
\hom{A}{\cohom{\coring{D}}{N}{X}}{Y}$, then $f$ is
$\coring{C}$--colinear if and only if the image of $(f \tensor{A}
N)\theta_X$ is included in $Y \cotensor{\coring{C}} N$. But these
are straightforward computations in view of the definition of the
contensor product. \\
2. Since the inclusion $\coring{C} \cotensor{\coring{C}} N
\subseteq \coring{C} \tensor{A} N$ splits-off, we get from
\ref{coass1} the natural isomorphism $(- \tensor{A} \coring{C})
\cotensor{\coring{C}} N \cong - \tensor{A} ( \coring{C}
\cotensor{\coring{C}} N) \cong - \tensor{A} N$. Now, the functor $
- \tensor{A} \coring{C}$ is right adjoint \cite[Lemma 3.1]
{Brzezinski:2000unp} to the forgetful functor $U_A :
\rcomod{\coring{C}} \rightarrow \rcomod{A}$ and, by hypothesis, $
- \cotensor{\coring{C}} N$ is right adjoint to the functor
$\cohom{\coring{D}}{N}{-} : \rcomod{\coring{D}} \rightarrow
\rcomod{\coring{C}}$. This implies that $- \tensor{A} N$ is right
adjoint to $U_A \cohom{\coring{D}}{N}{-}$, as desired.
\end{proof}
\end{punto}

\begin{example}
Let $N \in \bcomod{\coring{C}}{B}$. Then $N_B$ is quasi-finite if
and only if $- \tensor{A} N: \rcomod{A} \rightarrow \rcomod{B}$
has a left adjoint, that is, if and only if ${}_AN$ is finitely
generated and projective. In such a case, the left adjoint is $-
\tensor{B} \hom{A}{N}{A} : \rcomod{B} \rightarrow
\rcomod{\coring{C}}$. Notice that taking $B = K$, we obtain a
canonical structure of right $\coring{C}$--comodule on $N^* =
\hom{A}{N}{A}$ for every left $\coring{C}$--comodule $N$ such
that $N$ finitely generated and projective as a left $A$--module.
\end{example}

\begin{example}\label{restr}
Given an $A$--coring $\coring{C}$ and an $K$--algebra homomorphism
$\rho : A \rightarrow B$, we can consider the functor $-
\tensor{A} \coring{C}: \rcomod{B} \rightarrow \rcomod{\coring{C}}$
which is already the composite of the restriction of scalars
functor $(-)_{\rho} : \rcomod{B} \rightarrow \rcomod{A}$ and the
functor $- \tensor{A} \coring{C} : \rcomod{A} \rightarrow
\rcomod{\coring{C}}$. Since these functors have both left
adjoints, given, respectively, by the induction functor $ -
\tensor{A} B : \rcomod{A} \rightarrow \rcomod{B}$ and the
underlying functor $U_A : \rcomod{\coring{C}} \rightarrow
\rcomod{A}$, we get that the composite functor $- \tensor{A} B :
\rcomod{\coring{C}} \rightarrow \rcomod{B}$ is a left adjoint to
$- \tensor{A} \coring{C} : \rcomod{B} \rightarrow
\rcomod{\coring{C}}$. Clearly, $- \tensor{A} \coring{C} \cong -
\tensor{B} (B \tensor{A} \coring{C})$ and, thus, $B \tensor{A}
\coring{C}$ becomes a quasi-finite comodule.
\end{example}

\section{Separable homomorphisms of corings}\label{separable}

I propose a notion of homomorphism of corings which generalizes
both the concept of morphism of entwining structures
\cite{Brzezinski:1999unp} and the coring maps originally
considered in \cite{Sweedler:1975}. An induction functor is
constructed, which is shown to have a right adjoint, called
ad-induction functor. The separability of these two functors is
characterized in terms which generalize both the previous results
on rings \cite{Nastasescu/VandenBergh/VanOystaeyen:1989} and on
coalgebras \cite{Castano/Gomez/Nastasescu:1997}. Our approach
rests on the fundamental characterization of the separability of
adjoint functors given in \cite{Rafael:1990} and
\cite{DelRio:1992}.

Consider an $A$--coring $\coring{C}$ and a $B$--coring
$\coring{D}$, where $A$ and $B$ are $K$--algebras.

\begin{definition}
 A \emph{coring homomorphism} is a pair
$(\varphi,\rho)$, where $\rho: A \rightarrow B$ is a homomorphism
of $K$--algebras and $\varphi: \coring{C} \rightarrow \coring{D}$
is a homomorphism of $A$--bimodules and such that the following
diagrams are commutative:
\[
\xymatrix{
\coring{C} \ar^{\Delta_{\coring{C}}}[rr]
\ar^{\varphi}[dd] & & \coring{C}
\tensdia{A} \coring{C} \ar^{\varphi \tensfun{A} \varphi}[dr] & \\
 & & & \coring{D} \tensdia{A} \coring{D} \ar^{\omega_{\coring{D},\coring{D}}}[dl]\\
\coring{D} \ar^{\Delta_{\coring{D}}}[rr] & & \coring{D}
\tensdia{B} \coring{D} & } \qquad \xymatrix{ \coring{C}
\ar^{\epsilon_{\coring{C}}}[rr]
\ar^{\varphi}[dd]& & A \ar^{\rho}[dd] \\
& & & \\
\coring{D} \ar^{\epsilon_{\coring{D}}}[rr] & & B },
\]
where $\omega_{\coring{D},\coring{D}} : \coring{D} \tensdia{A}
\coring{D} \rightarrow \coring{D} \tensdia{B} \coring{D}$ is the
canonical map induced by $\rho : A \rightarrow B$.
\end{definition}

Throughout this section, we consider a coring homomorphism
$(\varphi,\rho) : \coring{C} \rightarrow \coring{D}$. We will
define the induction and ad-induction functors connecting the
categories of comodules $\rcomod{\coring{C}}$ and
$\rcomod{\coring{D}}$.

\begin{punto}\label{sigma}
We start with some unavoidable technical work. For every
$B$--bimodule $X$, let us denote by $\sigma_X : B \tensor{A} X
\rightarrow X \tensor{B} B$ the $B$--bimodule morphism given by $b
\tensor{A} x \mapsto bx \tensor{B} 1$. Given $B$--bimodules $X$,
$Y$, a straightforward computation shows that the diagram
\begin{equation}\label{xy1}
\xymatrix{ B \tensdia{A} X \tensdia{A} Y \ar^{B \tensfun{A}
\omega_{X,Y}}[d]
\ar^{\sigma_X \tensfun{A} Y}[rr]&  & X \tensdia{B} B \tensdia{A} Y \ar^{X \tensfun{B} \sigma_Y}[d]\\
B \tensdia{A} X \tensdia{B} Y \ar^{\sigma_{X \tensdia{B} Y}}[rr] &
& X \tensdia{B} Y \tensdia{B} B}
\end{equation}
commutes, where $\omega_{X,Y} : X \tensor{A} Y \rightarrow X
\tensor{B} Y$ is the obvious map. We have as well that for every
homomorphism of $B$--bimodules $f : X \rightarrow Y$ the following
diagram is commutative:
\begin{equation}\label{xy2}
\xymatrix{ B \tensdia{A} X \ar^{B \tensfun{A} f}[rr]
\ar^{\sigma_X}[d]& &
B \tensdia{A} Y \ar^{\sigma_Y}[d]\\
X \tensdia{B} B \ar^{f \tensfun{B} B }[rr] & & Y \tensdia{B} B}
\end{equation}
\end{punto}

\begin{punto}\label{labm}\textbf{The induction functor.}
Let $\lambda_M : M \rightarrow \coring{C} \tensor{A} M$ be a left
$\coring{C}$--comodule. Define $\widetilde{\lambda}_M : M
\rightarrow \coring{D} \tensor{A} M$ as the composite map
\[
\xymatrix{ M \ar_{\widetilde{\lambda}_M}[rr] \ar_{\lambda_M}[dr]&
& \coring{D}
\tensdia{A} M \\
& \coring{C} \tensdia{A} M \ar_{\varphi \tensfun{A} M}[ur]&}
\]
and $\lambda_{B \tensor{A} M} : B \tensor{A} M \rightarrow
\coring{D} \tensor{B} B \tensor{A} M$ as
\[
\xymatrix{ B \tensdia{A} M \ar^{\lambda_{B \tensdia{A} M}}[rr]
\ar_{B \tensfun{A} \widetilde{\lambda}_M}[dr]& & \coring{D} \tensdia{B} B \tensdia{A} M\\
 & B \tensdia{A} \coring{D} \tensdia{A} M \ar_{\sigma_{\coring{D}} \tensfun{A} M}[ur]& }
\]
\begin{proposition*}
The homomorphism of left $B$--modules $\lambda_{B \tensor{A} M}$
endows $B \tensor{A} M$ with a structure of left
$\coring{D}$--comodule. This gives a functor $B \tensor{A} - :
\lcomod{\coring{C}} \rightarrow \lcomod{\coring{D}}$.
\end{proposition*}
\begin{proof}
Consider the diagram
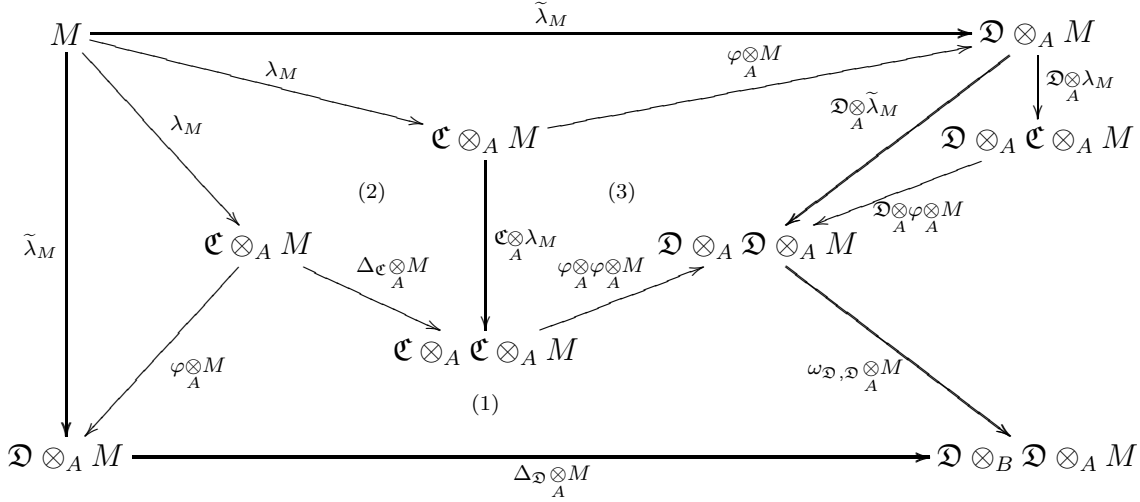
\[
\xymatrix{ M \ar^{\lambda_M}[ddr]
\ar^{\widetilde{\lambda}_{M}}[rrrr] \ar@<0.4pt>[rrrr]
\ar_{\widetilde{\lambda}_M}[dddd] \ar@<0.4pt>[dddd]
\ar^{\lambda_M}[drr]& \ar@{}[dddr]|{(2)} & \ar@{}[dddr]|{(3)}&
&\coring{D} \tensdia{A} M \ar_{\coring{D} \tensfun{A}
\widetilde{\lambda}_M}[ddl] \ar@<0.4pt>[ddl]
\ar^{\coring{D} \tensfun{A} \lambda_M}[d]\\
 & & \coring{C} \tensdia{A} M \ar^{\varphi \tensfun{A} M}[urr]
 \ar^{\coring{C} \tensfun{A} \lambda_M}[dd] & &
 \coring{D} \tensdia{A} \coring{C} \tensdia{A} M \ar^{\coring{D} \tensfun{A} \varphi \tensfun{A} M}[dl]\\
 &\coring{C} \tensdia{A} M \ar^{\varphi \tensfun{A} M}[ddl] \ar^{\Delta_{\coring{C}} \tensfun{A} M }[dr] &
 &
 \coring{D} \tensdia{A} \coring{D} \tensdia{A} M
 \ar_{\omega_{\coring{D},\coring{D}} \tensfun{A} M}[ddr]
 \ar@<0.4pt>[ddr]& \\
 & & \coring{C} \tensdia{A} \coring{C} \tensdia{A} M
 \ar^{\varphi \tensfun{A} \varphi \tensfun{A} M}[ur] \ar@{}[d]|{(1)} & & \\
\coring{D} \tensdia{A} M \ar_{\Delta_{\coring{D}} \tensfun{A} M
}[rrrr] \ar@<0.4pt>[rrrr]& & & & \coring{D} \tensdia{B} \coring{D}
\tensdia{A} M}
\]
The pentagon labeled as (1) commutes since $\varphi$ is a
homomorphism of corings. The four-edged diagram (2) is commutative
since $M$ is a left comodule. The commutation of the quadrilateral
(3) follows easily from the displayed decomposition of $\coring{D}
\tensor{A} \widetilde{\lambda}_M$. Therefore, the pentagon with
bold arrows commutes. Now consider the diagram
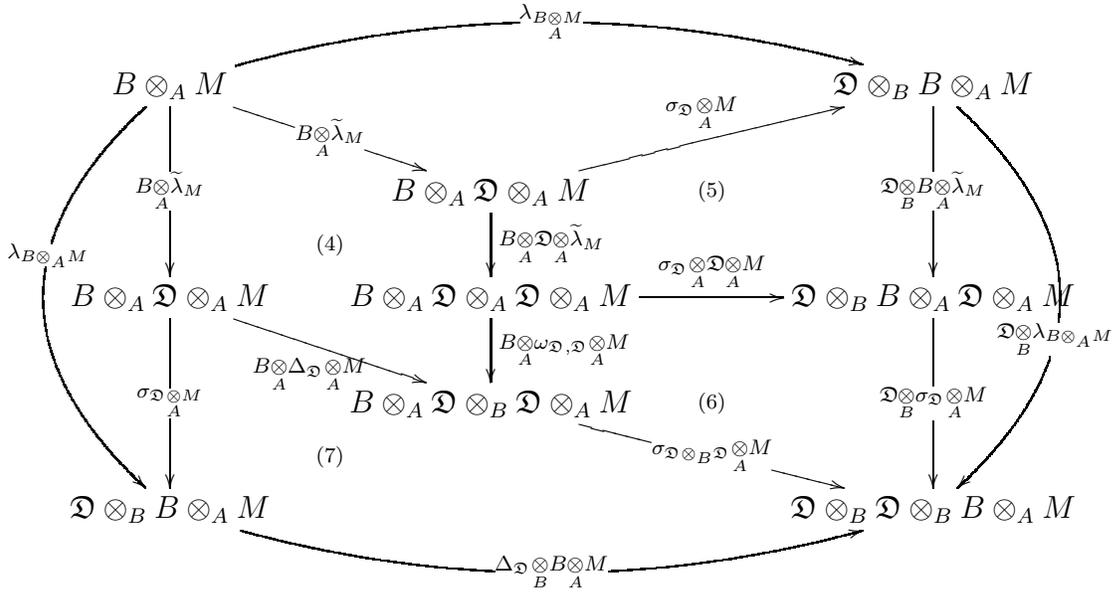
\[
\xymatrix{
B \tensor{A} M \ar@/^2pc/[rrr]|{\lambda_{B \tensfun{A} M}} 
\ar[dr]|{B \tensfun{A} \widetilde{\lambda}_M} 
\ar[dd]|{B \tensfun{A} \widetilde{\lambda}_M} 
\ar@/_4pc/[dddd]|(.4){\lambda_{B \tensor{A} M}} 
\ar@{}[dddr]|{(4)} &
 & &   \coring{D} \tensor{B} B \tensor{A} M
\ar[dd]|{\coring{D} \tensfun{B} B \tensfun{A} \widetilde{\lambda}_M} 
\ar@/^4pc/[dddd]|(.6){\coring{D} \tensfun{B} \lambda_{B \tensor{A} M}} \\ 
& B \tensor{A} \coring{D} \tensor{A} M
\ar^{\sigma_{\coring{D}} \tensfun{A} M}[urr] 
\ar^{B \tensfun{A} \coring{D} \tensfun{A}
\widetilde{\lambda}_M}[d] 
\ar@{}[rr]|{(5)}&
& &  \\ 
B \tensor{A} \coring{D} \tensor{A} M
\ar[dd]|{\sigma_{\coring{D} \tensfun{A} M}} 
\ar_{B \tensfun{A} \Delta_{\coring{D}} \tensfun{A} M}[dr] & B
\tensor{A} \coring{D} \tensor{A} \coring{D} \tensor{A} M
\ar^{B \tensfun{A} \omega_{\coring{D},\coring{D}} \tensfun{A} M}[d] 
\ar^{\sigma_{\coring{D}} \tensfun{A} \coring{D} \tensfun{A} M}[rr]
& & \coring{D} \tensor{B} B \tensor{A} \coring{D} \tensor{A} M
\ar[dd]|{\coring{D} \tensfun{B} \sigma_{\coring{D}} \tensfun{A} M}
  \\ 
& B \tensor{A} \coring{D} \tensor{B} \coring{D} \tensor{A} M
\ar[drr]|{\sigma_{\coring{D} \tensor{B} \coring{D}} \tensfun{A} M} 
\ar@{}[rr]|{(6)} & &  \\ 
\coring{D} \tensor{B} B \tensor{A} M
\ar@/_2pc/[rrr]|{\Delta_{\coring{D}}
\tensfun{B} B \tensfun{A} M} 
\ar@{}[ur]|{(7)} & & &   \coring{D} \tensor{B} \coring{D}
\tensor{B} B \tensor{A} M }
\]

We have proved before that (4) is commutative. Moreover, (5) is
obviously commutative and (6) and (7) commute by \eqref{xy1} and
\eqref{xy2}. It follows that the outer curved diagram commutes,
which gives the pseudo-coassocitative property for the coaction
$\lambda_{B \tensor{A} M}$. To check the counitary property, let
$\iota_M : M \rightarrow A \tensor{A} M$ be the canonical
isomorphism. We get from the commutative diagram
\[
\xymatrix{ M \ar^{\widetilde{\lambda}_M}[rr] \ar_{\iota_M}[dd]
\ar_{\lambda_M}[dr] & &
\coring{D} \tensor{A} M \ar^{\epsilon_{\coring{D}} \tensfun{A} M}[dd] \\
& \coring{C} \tensor{A} M \ar^{\epsilon_{\coring{C}} \tensfun{A}
M}[dl]
\ar_{\varphi \tensfun{A} M}[ur] & \\
A \tensor{A} M \ar_{\rho \tensfun{A} M}[rr] & & B \tensor{A} M }
\]
that the diagram
\[
\xymatrix{ B \tensor{A} M \ar^{\lambda_{B \tensor{A} M}}[rrrr] 
\ar|{B \tensfun{A} \widetilde{\lambda}_M}[drr] 
\ar@/_1pc/|(.7){\iota_{B \tensor{A} M}}[dddrrrr] 
\ar^{B \tensfun{A} \iota_M}[ddd] 
& & & & \coring{D} \tensor{B} B \tensor{A} M
\ar^{\epsilon_{\coring{D}} \tensfun{B} B \tensfun{A} M}[ddd] \\
& & B \tensor{A} \coring{D} \tensor{A} M \ar_{\sigma_{\coring{D}} \tensfun{A} M}[urr] 
\ar|(.6){B \tensfun{A} \epsilon_{\coring{D}} \tensfun{A} M}[dd] & & \\
& & & & \\
B \tensor{A} A \tensor{A} M \ar_{B \tensfun{A} \rho \tensfun{A}
M}[rr] & & B \tensor{A} B \tensor{A} M \ar_{\sigma_{B} \tensfun{A}
M}[rr] & & B \tensor{B} B \tensor{A} M}
\]
commutes, where $\iota_{B \tensor{A} M} : B \tensor{A} M
\rightarrow B \tensor{B} B \tensor{A} M$ denotes the canonical
isomorphism. Therefore, $\lambda_{B \tensor{A} M} : B \tensor{A} M
\rightarrow \coring{D} \tensor{B} B \tensor{A} M$ is a left
$\coring{D}$--comodule structure map. In order to show that the
assignment $M \mapsto B \tensor{A} M$ is functorial, we will prove
that $B \tensor{A} f$ is a homomorphism of left
$\coring{D}$--comodules for every morphism $f : M \rightarrow N$
in $\lcomod{\coring{C}}$. So, we have to show that the outer
rectangle in the following diagram is commutative
\[
\xymatrix{ B \tensor{A} M \ar^{B \tensfun{A} f}[rrrr] 
\ar|{B \tensfun{A} \widetilde{\lambda}_M}[dr] 
\ar|{\lambda_{B \tensor{A} M}}[dd]& & & &
B \tensor{A} N \ar|{B \tensor{A} \widetilde{\lambda}_N}[dl] 
\ar|{\lambda_{B \tensor{A} N}}[dd] \\
&  B \tensor{A} \coring{D} \tensor{A} M 
\ar|{\sigma_{\coring{D}} \tensfun{A} M}[dl] 
\ar^{B \tensfun{A} \coring{D} \tensfun{A} f}[rr]
& & B  \tensor{A} \coring{D} \tensor{A} N \ar|{\sigma_{\coring{D}} \tensfun{A} N}[dr]  & \\
\coring{D} \tensor{B} B \tensor{A} M 
\ar_{\coring{D} \tensfun{B} B \tensfun{A} f}[rrrr] 
 & & & & \coring{D} \tensor{B} B \tensor{M} }
\]
From the definition of $\widetilde{\lambda}_M,
\widetilde{\lambda}_N$ and the fact that $f$ is a morphism of
$\coring{C}$--comodules, it follows that the upper trapezium
commutes. The lower trapezium is commutative by \eqref{xy2}. Since
the two triangles commute by definition, we get that the outer
rectangle is commutative and $B \tensor{A} f$ is a morphism in
$\lcomod{\coring{D}}$.
\end{proof}
\end{punto}

\begin{punto}\label{numb}
Proposition \ref{labm} also implies by symmetry that for every
right $\coring{C}$--comodule $\rho_{M} : M \rightarrow M
\tensor{A} \coring{C}$, the right $B$--module $M \tensor{A} B$ is
endowed with a right $\coring{D}$--comodule structure $\rho_{ M
\tensor{A} B } : M \tensor{A} B \rightarrow  M \tensor{A} B
\tensor{B} \coring{D}$ given by $\rho_{M \tensor{A} B} = (M
\tensor{A} \delta_{\coring{D}})(\widetilde{\rho}_M \tensor{A} B)$,
where $\widetilde{\rho}_M = (M \tensor{A} \varphi) \rho_M$ and
$\delta_{\coring{D}} : \coring{D} \tensor{A} B \rightarrow B
\tensor{B} \coring{D}$ is the obvious map. We can already state
that
\begin{proposition*}
The assignment $M \mapsto M \tensor{A} B$ establishes a functor $-
\tensor{A} B : \rcomod{\coring{C}} \rightarrow
\rcomod{\coring{D}}$.
\end{proposition*}
\end{punto}

\begin{punto}\textbf{The ad-induction functor.}
Let us consider the left $\coring{D}$--comodule structure
\[
\xymatrix{\lambda_{B \tensor{A} \coring{C}} : B \tensor{A}
\coring{C} \ar[r] & \coring{D} \tensor{B} B \tensor{A} \coring{C}}
\]
defined on $B \tensor{A} \coring{C}$ in \ref{labm}.

We have as well a canonical structure of right
$\coring{C}$--comodule
\[
\xymatrix{ B \tensor{A} \Delta_{\coring{C}} : B \tensor{A}
\coring{C} \ar[r] & B \tensor{A} \coring{C} \tensor{A} \coring{C}
}
\]
\begin{proposition*}
The $B-A$--bimodule $B \tensor{A} \coring{C}$ is a
$\coring{D}-\coring{C}$--bicomodule, which is quasi-finite as a
right $\coring{C}$--comodule. Therefore, if ${}_A\coring{C}$
preserves the equalizer of $(\rho_Y \tensor{B} B \tensor{A}
\coring{C},Y \tensor{B} \lambda_{B \tensor{A} \coring{C}})$ for
every  right $\coring{D}$--comodule $Y$, then the functor
\[
\xymatrix{- \cotensor{\coring{D}} (B \tensor{A} \coring{C})
:\rcomod{\coring{D}} \ar[r] & \rcomod{\coring{C}}}
\]
is right adjoint to
\[
\xymatrix{- \tensor{A} B : \rcomod{\coring{C}} \ar[r] &
\rcomod{\coring{D}}}.
\]
\end{proposition*}
\begin{proof}
Since the comultiplication $\Delta_{\coring{C}}$ is coassociative,
we get that the diagram
\[
\xymatrix{\coring{C} \ar^{\widetilde{\lambda}_{\coring{C}}}[rr] 
\ar_{\Delta_{\coring{C}}}[dr] 
\ar^{\Delta_{\coring{C}}}[ddd] & & \coring{D} \tensor{A}
\coring{C} \ar^{\coring{D} \tensfun{A} \Delta_{\coring{C}}}[ddd]\\
& \coring{C} \tensor{A} \coring{C} \ar_{\varphi \tensfun{A} \coring{C}}[ur] %
\ar^{\coring{C} \tensfun{A} \Delta_{\coring{C}}}[d]& \\
& \coring{C} \tensor{A} \coring{C} \tensor{A} \coring{C} 
\ar^{\varphi  \tensfun{A} \coring{C} \tensfun{A} \coring{C}}[dr]& \\
\coring{C} \tensor{A} \coring{C} \ar^{\Delta_{\coring{C}} \tensfun{A} \coring{C}}[ur] 
\ar^{\widetilde{\lambda}_{\coring{C}} \tensfun{A} \coring{C}}[rr]&
& \coring{D} \tensor{A} \coring{C} \tensor{A} \coring{C}}
\]
is commutative. This implies that the left trapezium of the
following diagram is commutative.
\[
\xymatrix{ B \tensor{A} \coring{C} \ar^{\lambda_{B \tensor{A} \coring{C}}}[rr] 
\ar_{B \tensfun{A} \widetilde{\lambda}_{\coring{C}}}[dr] 
\ar^{B \tensfun{A} \Delta_{\coring{C}}}[ddd]
& & \coring{D} \tensor{B} B \tensor{A} \coring{C} %
\ar^{\coring{D} \tensfun{B} B \tensfun{A} \Delta_{\coring{C}}}[ddd] \\
& B \tensor{A} \coring{D} \tensor{A} \coring{C} %
\ar^{B \tensfun{A} \coring{D} \tensfun{A} \Delta_{\coring{C}}}[d] 
\ar_{\sigma_{\coring{D}} \tensfun{A} \coring{C}}[ur]&  \\
& B \tensor{A} \coring{D} \tensor{A} \coring{C} \tensor{A} \coring{C} %
\ar^{\sigma_{\coring{D}} \tensfun{A} \coring{C} \tensfun{A} \coring{C}}[dr]&  \\
B \tensor{A} \coring{C} \tensor{A} \coring{C} %
\ar^{B \tensfun{A} \widetilde{\lambda}_{\coring{C}} \tensfun{A} \coring{C}}[ur] 
\ar^{\lambda_{B \tensor{A} \coring{C}} \tensfun{A}
\coring{C}}[rr]& & \coring{D} \tensor{B} B \tensor{A} \coring{C}
\tensor{A} \coring{C}}
\]
Since we know that the rest of inner diagrams are also
commutative, we obtain that the outer diagram commutes, too. This
proves that $B \tensor{A} \coring{C}$ is a
$\coring{D}-\coring{C}$--bicomodule.

We will see that $B \tensor{A} \coring{C}$ is a quasi-finite
 right $\coring{C}$--comodule, i.e., that $\xymatrix{- \tensor{A}
B : \rcomod{\coring{C}} \ar[r] & \rcomod{B}}$ is left adjoint to
$\xymatrix{- \tensor{B} B \tensor{A} \coring{C} : \rcomod{B}
\ar[r] & \rcomod{\coring{C}}}$.
 Now, these functors fit in the
commutative diagrams
\[
\xymatrix{ \rcomod{A} \ar[rr]^{- \tensfun{A} \coring{C}}& &
\rcomod{\coring{C}}
\\ \rcomod{B} \ar[u]^{- \tensfun{B} B}
\ar[urr]_{- \tensfun{B} B \tensfun{A} \coring{C} }& & } \qquad
\xymatrix{\rcomod{\coring{C}}  \ar[drr]_{- \tensfun{A} B}
\ar[rr]^{U_A} & & \rcomod{A} \ar[d]^{- \tensfun{A} B}\\ & &
\rcomod{B}},
\]
where $U_A$ denotes the forgetful functor. Since $- \tensor{A} B$
is left adjoint to $ - \tensor{B} B$ and $U_A$ is left adjoint to
$ - \tensor{A} \coring{C}$ we get the desired adjunction. The
rest of the proposition follows from Proposition \ref{adjcohom}.
\end{proof}

\begin{remark*}
This proposition applies in the case that ${}_A\coring{C}$ is flat
or when $\coring{D}$ is a coseparable $B$--coring in the sense of
\cite{Guzman:1989}.
\end{remark*}
\end{punto}

\begin{punto}\textbf{The unit.}
Let $\overline{\Delta}_{\coring{C}} : \coring{C} \rightarrow
\coring{C} \tensor{A} B \tensor{B} B \tensor{A} \coring{C}$ be
the composite map
\[
\xymatrix{ \coring{C} \ar[r]^(.3){\Delta_{\coring{C}}} &
\coring{C} \tensor{A} \coring{C} \ar[r]^(.3){\iota} & \coring{C}
\tensor{A} B \tensor{B} B \tensor{A} \coring{C} }
\]
where $\iota$ maps $c \tensor{} c' \in \coring{C} \tensor{A}
\coring{C}$ to $c \tensor{} 1 \tensor{} 1 \tensor{} c'$. This map
is a homomorphism of $\coring{C}$--bicomodules. The unit $\theta$
of the adjunction $- \tensor{A} B \dashv - \tensor{B} B \tensor{A}
\coring{C}$ at a right $\coring{C}$--comodule $M$ is given by the
composite map
\[
\xymatrix{M \ar[r]^(.3){\rho_M} & M \tensor{A} \coring{C}
\ar[rr]^(.4){\varsigma_M \tensfun{A} \coring{C}} & & M \tensor{A}
B \tensor{B} B \tensor{A} \coring{C}}
\]
where $\varsigma_M : M \rightarrow M \tensor{A} B \tensor{B} B$
maps $m \in M$ to $m \tensor{A} 1 \tensor{B} 1$. We see, in
particular, that $\theta_{\coring{C}} =
\overline{\Delta}_{\coring{C}}$. By Proposition \ref{adjcohom},
$\theta_M$ factorizes throughout $(M \tensor{A} B )
\cotensor{\coring{D}} (B \tensor{A} \coring{C})$ and, therefore,
it gives the unit $\theta_M$ for the adjunction $ - \tensor{A} B
\dashv - \cotensor{\coring{D}} ( B \tensor{A} \coring{C})$ at $M$.
So, the multiplication $\Delta_{\coring{C}}$ finally induces a map
\[
\xymatrix{\overline{\Delta}_{\coring{C}} : \coring{C} \ar[r] &
(\coring{C} \tensor{A} B) \cotensor{\coring{D}} (B \tensor{A}
\coring{C})}
\]
which is a homomorphism of $\coring{C}$--bicomodules.
\end{punto}

We are now ready to state the characterization of the separability
of the induction functor.

\begin{theorem}\label{sepind}
Assume that $\coring{A}_A$ preserves the equalizer of $(\rho_Y
\tensor{B} B \tensor{B} \tensor{A} \coring{C}, Y \tensor{B}
\lambda_{B \tensor{A} \coring{C}})$ for every $Y \in
\rcomod{\coring{D}}$, and that $X_A$ preserves the equalizer of
$(\rho_{\coring{C} \tensor{A} B} \tensor{B} B \tensor{A}
\coring{C}, \coring{C} \tensor{A} B \tensor{B} \lambda_{B
\tensor{A} \coring{C}})$ for every $X \in \rcomod{\coring{C}}$.
The functor $ \xymatrix{- \tensor{A} B : \rcomod{\coring{C}}
\ar[r] & \rcomod{\coring{D}}}$ is separable if and only if there
is a homomorphism if $\coring{C}$--bicomodules
\[
\xymatrix{\omega_{\coring{C}} : (\coring{C} \tensor{A} B)
\cotensor{\coring{D}} (B \tensor{A} \coring{C}) \ar[r] &
\coring{C}}
\]
such that $\omega_{\coring{C}} \overline{\Delta}_{\coring{C}} =
\coring{C}$.
\end{theorem}
\begin{proof}
Assume that $ - \tensor{A} B$ is separable. By \cite[Theorem
4.1]{DelRio:1992} or \cite[Theorem 1.2]{Rafael:1990}, the unit of
the adjunction
\[
\xymatrix{\theta : 1_{\rcomod{\coring{C}}} \ar[r] & (- \tensor{A}
B) \cotensor{\coring{D}} (B \tensor{A} \coring{D})}
\]
is split-mono, that is, there is a natural transformation
$\xymatrix{\omega : (- \tensor{A} B) \cotensor{\coring{D}} (B
\tensor{A} \coring{C}) \ar[r]
 & 1_{\rcomod{\coring{C}}}}$ such that $\omega \theta =
1_{\rcomod{\coring{C}}}$. By Lemma \ref{coass1} the functor
$(-\tensor{A} B) \cotensor{\coring{D}} (B \tensor{A} \coring{C})$
is $\coring{C}$--compatible. By Proposition \ref{hacebicom},
$\omega_{\coring{C}}$ is a $\coring{C}$--bicomodule map.
Obviously, $\omega_{\coring{C}} \overline{\Delta}_{\coring{C}} =
\coring{C}$.\\
To prove the converse, we need to construct a natural
transformation $\omega$ from the bicomodule map
$\omega_{\coring{C}}$. Given a right $\coring{C}$--comodule $M$,
consider the diagram
\[
\xymatrix{ M \tensor{A} B \tensor{B} \coring{D} \tensor{B} B \tensor{A} \coring{C}
\ar^{\rho_M \tensfun{A} B \tensfun{B} \coring{D} \tensfun{B} B
\tensfun{A} \coring{C}}[rrrr]&
 & &
  & M \tensor{A} \coring{C} \tensor{A} B \tensor{B} \coring{D} \tensor{B} B \tensor{A} \coring{C}\\
   & & & & \\
M \tensor{A} B \tensor{B} B \tensor{A} \coring{C}
\ar^{\rho_M \tensfun{A} B \tensfun{B} B \tensfun{A}
\coring{C}}[rrrr]
\ar^{\rho_{M \tensor{A} B} \tensfun{B} B \tensfun{A}
\coring{C}}@<1ex>[uu]
\ar_{M \tensfun{A} B \tensfun{B} \lambda_{B \tensor{A}
\coring{C}}}@<-1ex>[uu]&
& & &  M \tensor{A} \coring{C} \tensor{A} B \tensor{B} B \tensor{A} \coring{C}
\ar^{M \tensfun{A} \rho_{\coring{C} \tensor{A} B} \tensfun{B} B \tensfun{A} \coring{C}}@<1ex>[uu]
\ar_{M \tensor{A} \coring{C} \tensor{A} B \tensor{B} \lambda_{B \tensor{A} \coring{C}}}@<-1ex>[uu]\\
 & & & & \\
(M \tensor{A} B) \cotensor{\coring{D}} (B \tensor{A} \coring{C})
\ar@{..>}^{\kappa_M}[rrrr] \ar[uu]& & & & M \tensor{A}
((\coring{C} \tensor{A} B) \cotensor{\coring{D}} (B \tensor{A}
\coring{C})) \ar[uu]},
\]
where the vertical are equalizer diagrams (here, we are using the
definition of the cotensor product and the fact that $M_A$
preserves the equalizer of $(\rho_{\coring{C} \tensor{A} B}
\tensor{B} B \tensor{A} \coring{C}, \coring{C} \tensor{A} B
\tensor{B} \lambda_{B \tensor{A} \coring{C}})$). If we prove that
the top rectangle commutes, then there is a unique dotted arrow
$\kappa_M$ making the bottom rectangle commute. The identity
\[
(M \tensfun{A} \coring{C} \tensfun{A} B \tensfun{B} \lambda_{B
\tensor{A} \coring{C}})(\rho_M \tensfun{A} B \tensfun{B} B
\tensfun{A} \coring{C}) = (\rho_M \tensfun{A} B \tensfun{B}
\coring{D} \tensfun{B} B \tensfun{A} \coring{C})(M \tensfun{A} B
\tensfun{B} \lambda_{B \tensor{A }\coring{C}})
\]
is obvious; so we need just to prove that
\[
(M \tensor{A} \rho_{\coring{C} \tensor{A} B} \tensfun{B} B
\tensfun{A} \coring{C})(\rho_M \tensfun{A} B \tensfun{B} B
\tensfun{A} \coring{C}) = (\rho_M \tensfun{A} B \tensfun{B} B
\tensfun{B} \coring{D} \tensfun{B} B \tensfun{A}
\coring{C})(\rho_{M \tensor{A} B} \tensfun{B} B \tensfun{A}
\coring{C}),
\]
which is equivalent to
\[
(M \tensfun{A} \rho_{\coring{C} \tensor{A} B})(\rho_M \tensfun{A}
B) = ( \rho_M \tensfun{A} B \tensfun{B} \coring{D})\rho_{M
\tensor{A} B},
\]
and this last identity is easy to check. Now, consider the diagram
\[
\xymatrix{ M \ar[ddd]^{\rho_M} \ar@{=}[dr] \ar[rrr]^{\theta_M} & &
&
 M \tensor{A} B \tensor{B} B \tensor{A} \coring{C} 
\ar[ddd]_{\rho_M \tensfun{A} B \tensfun{B} B \tensfun{A}
\coring{C}}\\
& M \ar[d]^{\rho_M} \ar[r]^(.3){\theta_M} & (M \tensor{A} B)
\cotensor{\coring{D}} (B \tensor{A}
\coring{C}) \ar@{^{(}->}[ur] \ar[d]^{\kappa_M} &  \\
& M \tensor{A} \coring{C} \ar[r]^(.3){M \tensor{A}
\overline{\Delta}_{\coring{C}}} \ar@{=}[dl] & M \tensor{A}
((\coring{C} \tensor{A} B) \cotensor{\coring{D}} (B \tensor{A}
\coring{C}))
\ar@{^{(}->}[dr] & \\
M \tensor{A} \coring{C} \ar[rrr]^{M \tensfun{A}
\overline{\Delta}_{\coring{C}}} & & & M \tensor{A} \coring{C}
\tensor{A} B \tensor{B} B \tensor{A} \coring{C}}
\]

By Lemma \ref{FunoFdos}, we have that $M \tensor{A} \theta_{M} =
\theta_{M \tensor{A} \coring{C}}$. Since $\theta$ is natural, and
$\theta_M = \overline{\Delta}_{\coring{C}}$, this implies that the
external rectangle commutes. We know that the four trapezia
commute, whence the internal rectangle commutes as well. Define
$\omega_M = (M \tensor{A} \epsilon_{\coring{C}})(M \tensor{B}
\omega_{\coring{C}})\kappa_M$, which gives a natural
transformation $\omega : (- \tensor{A} B) \cotensor{\coring{D}} (B
\tensor{A} \coring{C})  \rightarrow 1_{\rcomod{\coring{C}}}$.
Moreover,
\begin{multline*}
\omega_M \theta_M = (M \tensor{A} \epsilon_{\coring{C}})(M
\tensor{A} \omega_{\coring{C}})\kappa_M\theta_M = \\ (M \tensor{A}
\epsilon_{\coring{C}})(M \tensor{A}\omega_{\coring{C}})(M
\tensor{A} \theta_M)\rho_M = (M \tensor{A} \epsilon_M)\rho_M = M
\end{multline*}
Therefore, $\omega \theta = 1_{\rcomod{\coring{C}}}$ and, by
\cite[Theorem 1.2]{Rafael:1990}, the functor $- \tensor{A} B$ is
separable.
\end{proof}

The counit map $\epsilon_{\coring{C}} : \coring{C} \rightarrow A$
is a homomorphism of $A$--corings, where we consider on $A$ the
canonical $A$--coring structure. When applied to
$\epsilon_{\coring{C}}$, Theorem \ref{sepind} boils down to

\begin{corollary}\cite[Theorem 3.5]{Brzezinski:2000unp}\label{sepolvido}
For an $A$--coring $\coring{C}$, the forgetful functor $U_A :
\rcomod{\coring{C}} \rightarrow \rcomod{A}$ is separable if and
only if there is an $A$--bimodule map $\gamma : \coring{C}
\tensor{A} \coring{C} \rightarrow A$ such that $\gamma
\Delta_{\coring{C}} = \epsilon_{\coring{C}}$ and, in Sweedler's
sigma notation,
\[
c_{(1)}\gamma(c_{(2)} \tensor{A} c') = \gamma  (c \tensor{A}
c'_{(1)})c'_{(2)} \qquad \text{for all } \; c, c' \in \coring{C}
\]
\end{corollary}
\begin{proof}
Obviously, the forgetful functor coincides with $ - \tensor{A} A$,
so that we get from Theorem \ref{sepind} a characterization of the
separability of this functor in terms of the existence of a
$\coring{C}$--bicomodule map $\omega_{\coring{C}} : \coring{C}
\tensor{A} \coring{C} \rightarrow \coring{C}$ such that
$\omega_{\coring{C}} \Delta_{\coring{C}} = \coring{C}$. Now,
notice that the adjointness isomorphism
\[
\hom{\coring{C}}{\coring{C} \tensor{A} \coring{C}}{\coring{C}}
\cong \hom{A}{\coring{C} \tensor{A} \coring{C}}{A}
\]
transfers faithfully the mentioned properties of
$\omega_{\coring{C}}$ to the desired properties of $\gamma =
\epsilon_{\coring{C}}\omega_{\coring{C}}$.
\end{proof}

\begin{punto}\textbf{The counit of the adjunction.}
Let $\hat{\epsilon}_{\coring{C}}$ be the homomorphism of
$B$--bimodules that makes commute the following diagram
\[
\xymatrix{B \tensor{A} \coring{C} \tensor{A} B
\ar[rr]^{\hat{\epsilon}_{\coring{C}}} 
\ar[d]^{B \tensor{A} \epsilon_{\coring{C}} \tensor{A} B} & & B \\
B \tensor{A} A \tensor{A} B \ar[rr]^{B \tensfun{A} \rho
\tensfun{A} B} & & B \tensor{A} B \tensor{A} B \ar[u]^{m}},
\]
where $m : B \tensor{A} B \tensor{A} B \rightarrow B$ is the
obvious multiplication map. Define, for every right $B$--module
$Y$, $\chi_Y = \mu_Y (Y \tensor{A} \hat{\epsilon}_{\coring{C}})$,
where $\mu_Y : Y \tensor{B} B \rightarrow Y$ is the canonical
isomorphism. This natural transformation $\chi$ gives the counit
of the adjunction $- \tensor{A} B \dashv - \tensor{B} ( B
\tensor{A} \coring{C})$. By \ref{adjcohom}, the counit of the
adjunction $- \tensor{A} B \dashv - \cotensor{\coring{D}} ( B
\tensor{A} \coring{C})$ is given by the restriction of $\chi$ to
$(- \cotensor{\coring{D}} (B \tensor{A} \coring{C}) \tensor{A} B$.
We shall use the same notation for this counit. Now, define
$\hat{\varphi}$ as the $B$--bimodule map completing the diagram
\[
\xymatrix{B \tensor{A} \coring{C} \tensor{A} B
\ar[r]^(.7){\hat{\varphi}} \ar[d]_{B \tensfun{A} \varphi \tensfun{A} B} & \coring{D} \\
B \tensor{A} \coring{D} \tensor{A} B \ar[ur]_{m_{\coring{D}}} & }
\]
where $m_{\coring{D}} : B \tensor{A} \coring{D} \tensor{A} B
\rightarrow \coring{D}$ is the obvious multiplication map given
by the $B$--bimodule structure of $\coring{D}$. We claim that
$\hat{\varphi}$ is a $\coring{D}$--comodule map. To prove this, we
first show that the diagram
\begin{equation}\label{phihat}
\xymatrix{(\coring{D} \cotensor{\coring{D}} (B \tensor{A}
\coring{C})) \tensor{A} B \ar[r] & \coring{D} \tensor{B} B
\tensor{A} \coring{C} \tensor{A} B \ar[d]^{\chi_{\coring{D}}} \\
B \tensor{A} \coring{C} \tensor{A} B \ar[u]^{\lambda_{B \tensor{A}
\coring{C}} \tensfun{A} B} \ar[r]^(.6){\hat{\varphi}} &
\coring{D}}
\end{equation}
is commutative. This is done by the following computation, where
$\widetilde{m} : B \tensor{B} B \tensor{A} B \rightarrow B$
denotes the obvious multiplication map:

\begin{multline*}
\chi_{\coring{D}} (\lambda_{B \tensor{A} \coring{C}} \tensfun{A} B
)  = \mu_{\coring{D}}(\coring{D} \tensfun{A}
\hat{\epsilon}_{\coring{C}})(\sigma_{\coring{D}} \tensfun{A}
\coring{C} \tensfun{A} B)(B \tensfun{A}
\widetilde{\lambda}_{\coring{C}} \tensfun{A} B) = \\
   \mu_{\coring{D}}(\coring{D} \tensfun{B} m)(\coring{D}
 \tensfun{B} B \tensfun{A} \rho \tensfun{A} B)(\coring{D} \tensfun{B} B
 \tensfun{A} \epsilon_{\coring{C}} \tensfun{A}
 \coring{D})(\sigma_{\coring{D}} \tensfun{A} \coring{C}
 \tensfun{A} B)(B \tensfun{A} \varphi \tensfun{A} \coring{C}
 \tensfun{A} B)(B \tensfun{A} \Delta_{\coring{C}} \tensfun{A} B)
 =\\
 \mu_{\coring{D}}(\coring{D} \tensfun{B} m)(\coring{D} \tensfun{B}
 B \tensfun{A} \epsilon_{\coring{D}} \tensfun{A} B)(\coring{D}
 \tensfun{B} B \tensfun{A} \varphi \tensfun{A}
 B)(\sigma_{\coring{D}} \tensfun{A} \coring{C} \tensfun{A} B)(B
 \tensfun{A} \varphi \tensfun{A} \coring{C} \tensfun{A} B)(B
 \tensfun{A} \Delta_{\coring{C}} \tensfun{A} B) = \\
\mu_{\coring{D}}(\coring{D} \tensfun{B} m)(\coring{D} \tensfun{B}
 B \tensfun{A} \epsilon_{\coring{D}} \tensfun{A}
 B)(\sigma_{\coring{D}} \tensfun{A} \varphi \tensfun{A} B)(B
 \tensfun{A} \varphi \tensfun{A} \coring{C} \tensfun{A} B)(B
 \tensfun{A} \Delta_{\coring{C}} \tensfun{A} B) = \\
 \mu_{\coring{D}}(\coring{D} \tensfun{B} m)(\coring{D} \tensfun{B}
 B \tensfun{A} \epsilon_{\coring{D}} \tensfun{A}
 B)(\sigma_{\coring{D}} \tensfun{A} \coring{D} \tensfun{A} B)(B
 \tensfun{A} \varphi \tensfun{A} \varphi \tensfun{A} B)(B
 \tensfun{A} \Delta_{\coring{C}} \tensfun{A} B)) = \\
 \mu_{\coring{D}}(\coring{D} \tensfun{B} \widetilde{m})(\coring{D}
 \tensfun{B} B \tensfun{B} \epsilon_{\coring{D}} \tensfun{A}
 B)(\sigma_{\coring{D}} \tensfun{B} \coring{D} \tensfun{A} B)(B
 \tensfun{A} \omega_{\coring{D},\coring{D}} \tensfun{A} B)(B
 \tensfun{A} \varphi \tensfun{A} \varphi \tensfun{A} B)(B
 \tensfun{A} \Delta_{\coring{D}} \tensfun{A} B) = \\
 \mu_{\coring{D}}(\coring{D} \tensfun{B} \widetilde{m})(\coring{D}
 \tensfun{B} B \tensfun{B} \epsilon_{\coring{D}} \tensfun{A}
 B)(\sigma_{\coring{D}} \tensfun{B} \coring{D} \tensfun{A} B)(B
 \tensfun{A} \Delta_{\coring{D}} \tensfun{A} B)(B \tensfun{A}
 \varphi \tensfun{A} B) = \\
 \mu_{\coring{D}}(\coring{D} \tensfun{B}
 \widetilde{m})(\sigma_{\coring{D}} \tensfun{B}
 \epsilon_{\coring{D}} \tensfun{A} B)(B \tensfun{A}
 \Delta_{\coring{D}} \tensfun{A} B)(B \tensfun{A} \varphi
 \tensfun{A} B) = \\
 \mu_{\coring{D}}(\coring{D} \tensfun{B} \widetilde{m})(\sigma_{\coring{D}} \tensfun{B} B \tensfun{A}
 B)(B \tensfun{A} \coring{D} \tensfun{B} \epsilon_{\coring{D}}
 \tensfun{A} B) ( B \tensfun{A} \Delta_{\coring{D}} \tensfun{A}
 B)(B \tensfun{A} \varphi \tensfun{A} B) = \\
 \mu_{\coring{D}}(\coring{D} \tensfun{B}
 \widetilde{m})(\sigma_{\coring{D}} \tensfun{B} B \tensfun{A} B)(B
 \tensfun{A} \varphi \tensfun{A} B) = \\
 \mu_{\coring{D}}(\coring{D} \tensfun{B} m)(\sigma_{\coring{D}} \tensfun{A} B \tensfun{A} B)(B
 \tensfun{A} \varphi \tensfun{A} B) = \hat{\varphi}
\end{multline*}
We know that $\lambda_{B \tensor{A} \coring{C}} \tensor{A} B$ is a
homomorphism of $\coring{D}$--bicomodules and, by \ref{hacebicom},
$\chi_{\coring{D}}$ is $\coring{D}$--bicolinear, too. This proves
that $\hat{\varphi} : B \tensor{A} \coring{C} \tensor{A} B
\rightarrow \coring{D}$ is a homomorphism of
$\coring{D}$--bicomodules.
\end{punto}

We are now in a position to prove our separability theorem for the
ad-induction functor.

\begin{theorem}\label{sepcores}
Assume that ${}_AB$ and ${}_A\coring{C}$ preserve the equalizer of
$(\rho_M \tensor{B} B \tensor{A} \coring{C}, M \tensor{B}
\lambda_{B \tensor{A} \coring{C}})$ for every right
$\coring{D}$--comodule $M$ (e.g. ${}_AB$ and ${}_A\coring{C}$ are
flat or $\coring{D}$ is a coseparable $B$--coring in the sense of
\cite{Guzman:1989}). The functor $\xymatrix{-
\cotensor{\coring{D}} (B \tensor{A} \coring{C}) :
\rcomod{\coring{D}} \ar[r] & \rcomod{\coring{C}}}$ is separable if
and only if there exists a $\coring{D}$--bicomodule homomorphism
$\xymatrix{\hat{\nu}_{\coring{D}} : \coring{D} \ar[r]& B
\tensor{A} \coring{C} \tensor{A} B }$ such that $\hat{\varphi}
\hat{\nu}_{\coring{D}} = \coring{D}$.
\end{theorem}
\begin{proof}
If $- \cotensor{\coring{D}} (B \tensor{A} \coring{C})$ is
separable then, by \cite[Theorem 1.2]{Rafael:1990}, there exists a
natural transformation
\[
\xymatrix{\nu : 1_{\rcomod{\coring{D}}} \ar[r] & ( -
\cotensor{\coring{D}} (B \tensor{A} \coring{C})) \tensor{A} B}
\]
such that $\chi \nu = 1_{\rcomod{\coring{D}}}$. In particular,
$\chi_{\coring{D}}\nu_{\coring{D}} = \coring{D}$ and, by
Proposition \ref{hacebicom}, $\nu_{\coring{D}}$ is a bicomodule
map (in fact, we easily get that the functor $(-
\cotensor{\coring{D}} (B \tensor{A} \coring{C})) \tensor{A} B$ is
$\coring{D}$--compatible from the fact that $\coring{D}
\cotensor{\coring{D}} (B \tensor{A} \coring{C})$ is a direct
summand of $\coring{D} \tensor{B} B \tensor{A} \coring{C}$ as a
left $B$--module). The map $\lambda_{B \tensor{A} \coring{C}}
\tensor{A} B$ gives an isomorphism of $\coring{D}$--bicomodules
$B \tensor{A} \coring{C} \tensor{A} B \cong (\coring{D}
\cotensor{D} (B \tensor{A} \coring{C})) \tensor{A} B$ which
implies, after \eqref{phihat}, that
\[
\coring{D} = \chi_{\coring{D}}\nu_{\coring{D}}  =
\hat{\varphi}(\lambda_{B \tensor{A} \coring{C}} \tensor{A}
B)^{-1}\nu_{\coring{D}}
\]
Thus, $\hat{\nu}_\coring{D} = (\lambda_{B \tensor{A} \coring{C}}
\tensor{A} B)^{-1}\nu_{\coring{D}}$ is the desired
$\coring{D}$--bicomodule map.

For the converse, assume there is a $\coring{D}$--bicomodule map
$\hat{\nu}_{\coring{D}} : \coring{D} \rightarrow B \tensor{A}
\coring{C} \tensor{A} B$ such that $\hat{\varphi}
\hat{\nu}_{\coring{D}} = \coring{D}$. For each right
$\coring{D}$--comodule $M$, let us prove that the homomorphism of
right $\coring{D}$--comodules
\begin{equation}\label{prenu}
\xymatrix{M \ar[r]^(.4){\rho_M} & M \tensor{B} \coring{D}
\ar[rr]^(.4){M \tensfun{B} \hat{\nu}_{\coring{D}}} & & M
\tensor{B} B \tensor{A} \coring{C} \tensor{A} B}
\end{equation}
factorizes throughout $(M \cotensor{\coring{D}} (B \tensor{A}
\coring{C})) \tensor{A} B$. Since ${}_AB$ preserves the equalizer
of $(\rho_M \tensor{B} B \tensor{A} \coring{C}, M \tensor{B}
\lambda_{B \tensor{A} \coring{C}})$, we know that
\[
(M \cotensor{\coring{D}} (B \tensor{A} \coring{C})) \tensor{A} B
\cong M \cotensor{\coring{D}} (B \tensor{A} \coring{C} \tensor{A}
B).
\]
Therefore, we need just to check the equality
\begin{equation}\label{Rolanda}
(\rho_M \tensor{B} B \tensor{A} \coring{C} \tensor{A} B)(M
\tensor{B} \hat{\nu}_{\coring{D}})\rho_M = (M \tensor{B}
\lambda_{B \tensor{A} \coring{C} \tensor{A} B})(M \tensor{B}
\hat{\nu}_{\coring{D}})\rho_M
\end{equation}
This is done by the following computation:
\begin{multline*}
(M \tensor{B} \lambda_{B \tensor{A} \coring{C} \tensor{A} B})(M
\tensor{B} \hat{\nu}_{\coring{D}})\rho_M = \\
(M \tensor{B} \coring{D} \tensor{B} \hat{\nu}_{\coring{D}})(M
\tensor{B} \Delta_{\coring{D}})\rho_M = (M \tensor{B} \coring{D}
\tensor{B}
\hat{\nu}_{\coring{D}})\rho_M = \\
(\rho_M \tensor{B} \hat{\nu}_{\coring{D}})\rho_M = (\rho_M
\tensor{B} B \tensor{A} \coring{C} \tensor{A} B)(M \tensor{B}
\hat{\nu}_{\coring{D}})\rho_M
\end{multline*}
Thus, we have proved that the natural transformation given in
\eqref{prenu} factorizes throughout a natural transformation
\begin{equation*}\label{nu}
\xymatrix{\nu_M : M \ar[r] & (M \cotensor{\coring{D}} (B
\tensor{A} \coring{C})) \tensor{A} B}
\end{equation*}
This means that we have a commutative diagram
\begin{equation}\label{rombo}
\xymatrix{ M \ar^{\rho_M}[r] \ar^{\nu_M}[drr]&
 M \tensor{B} \coring{D} \ar^(.4){M \tensfun{B} \hat{\nu}_{\coring{D}}}[r]&
  M \tensor{B} B \tensor{A} \coring{C} \tensor{A} B \ar^{\chi_M}[drr]& & \\
& & (M \cotensor{\coring{D}} (B \tensor{A} \coring{C}))
\tensor{A} B \ar_(.7){\chi_M}[rr] \ar[u]& &M }
\end{equation}
Finally, we will show that $\nu_M$ splits off $\chi_M$ by means
of the following computation:
\begin{alignat*}{2}
\rho_M \chi_M \nu_M & = \rho_M\chi_M(M \tensfun{B}
\hat{\nu}_{\coring{D}})\rho_M & \\
& = \chi_{M \tensor{B} \coring{D}}(\rho_M \tensfun{B} B
\tensfun{A} \coring{C} \tensfun{A} B)(M \tensfun{B}
\hat{\nu}_{\coring{D}})\rho_M & \qquad (\chi \; \text{is
natural})\\
 & = (M \tensfun{B} \chi_{\coring{D}})(\rho_M \tensfun{B} B
\tensfun{A} \coring{C} \tensfun{A} B)(M \tensfun{B}
\hat{\nu}_{\coring{D}})\rho_M & \qquad (\text{by Lemma
\ref{FunoFdos})} \\
 & = (M \tensfun{B} \chi_{\coring{D}})(M \tensfun{B} \lambda_{B
 \tensor{A} \coring{C} \tensor{A} B})(M \tensfun{B}
\hat{\nu}_{\coring{D}})\rho_M & \qquad (\hat{\nu}_{\coring{D}} \;
\text{is bicolinear}) \\
& = (M \tensfun{B} \hat{\varphi})(M \tensfun{B}
\hat{\nu}_{\coring{D}})\rho_M  = \rho_M &
\end{alignat*}

Since $\rho_M$ is a monomorphism, we get $\chi_M \nu_M = M$. By
\cite[Theorem 1.2]{Rafael:1990}, $- \cotensor{D} ( B \tensor{A}
\coring{C})$ is a separable functor.
\end{proof}

By applying the stated theorem to $\epsilon_{\coring{C}} :
\coring{C} \rightarrow A$ we obtain:

\begin{corollary}\cite[Theorem 3.3]{Brzezinski:2000unp}
Let $\coring{C}$ be an $A$--coring. Then the functor $- \tensor{A}
\coring{C}$ is separable if and only if there exists an invariant
$e \in \coring{C}$ (that is, $e \in \coring{C}$ satisfying $ae =
ea$ for every $a \in A$) such that $\epsilon_{\coring{C}}(e) = 1$.
\end{corollary}
\begin{proof}
This follows from Theorem \ref{sepcores} taking that the
$A$--bimodule homomorphisms from $A$ to $\coring{C}$ correspond
bijectively with the invariants of $\coring{C}$ into account.
\end{proof}
\begin{punto}\textbf{Final Remarks.}
T. Brzezi\'nski showed \cite[Proposition 2.2]{Brzezinski:2000unp}
that if $(A,C)_{\psi}$ is an entwining structure over $K$ then $A
\tensor{K} C$ can be endowed with an $A$--coring structure in
such a way that the category $\rcomod{A \tensor{K} C}$ is
isomorphic with the category $\cat{M}_A^C(\psi)$ of entwined
$(A,C)_{\psi}$--modules. It turns out that every morphism of
entwining structures $(f,g) : (A,C)_{\psi} \rightarrow
(B,D)_{\gamma}$ in the sense of \cite{Brzezinski:1999unp} gives a
coring morphism $(f \tensor{K} g, f) : A \tensor{K} C \rightarrow
B \tensor{K} D$. Some straightforward computations show that the
statements of the separability theorem \cite[Theorem
3.4]{Brzezinski:1999unp} correspond to our theorems \ref{sepind}
and \ref{sepcores}. In fact, one can give the notions of totally
integrable and totally cointegrable morphism of entwining
structures to the framework of morphisms of corings $(\phi,\rho):
\coring{C} \rightarrow \coring{D}$ by requiring the existence of
the splitting bicomodule maps $\omega_{\coring{C}}$ and
$\hat{\nu}_{\coring{D}}$, respectively. In the first case, the
existence of $\omega_{\coring{C}}$ can be transferred, if desired,
to the existence of certain bimodule map with extra properties by
means of the adjointness isomorphism
\[
\hom{\coring{C}}{(\coring{C} \tensor{A} B) \cotensor{\coring{D}}
(B \tensor{A} \coring{C})}{\coring{C}} \cong \hom{A}{(\coring{C}
\tensor{A} B) \cotensor{\coring{D}} (B \tensor{A} \coring{C})}{A}
\]
\end{punto}

\noindent \textbf{Acknowledgement:} I would like to thank Tomasz
Brzezi{\'n}ski for some helpful comments.

\providecommand{\bysame}{\leavevmode\hbox
to3em{\hrulefill}\thinspace}


\end{document}